\RequirePackage{fix-cm}
\documentclass[twocolumn]{svjour3}          
\usepackage[english]{babel}
\usepackage[cp1251]{inputenc}
\usepackage{amsmath,amssymb}
\usepackage{float}
\usepackage[framed]{mcode}         
\smartqed  
\usepackage{graphicx}
\usepackage[suffix=, outdir=./]{epstopdf}
\begin{document}

\title{A dynamical model of Sayano-Shushenskaya hydropower plant: stability, oscillations, and accident}
\subtitle{Dynamical model of Sayano-Shushenskaya hydropower plant}


\author{
	Leonov G.A.      \and
	Kuznetsov N.V.   \and
	Solovyeva E.P.
}


\institute{
Leonov G.A., Kuznetsov N.V., Solovyeva E.P. \at
Saint-Petersburg State University;
University of Jyvaskyla.\\
\email{nkuznetsov239@gmail.com (corr. author email)}
}


\date{Received: date / Accepted: date}

\maketitle

\begin{abstract}
This work is devoted to the construction and study of a mathematical model of hydropower unit,
consisting of synchronous generator, hydraulic turbine, and speed governor.
It is motivated by the accident happened on the Sayano-Shushenskaya
hydropower plant in 2009 year. Parameters of the Sayano-Shushenskaya hydropower plant
were used for modeling the system. Oscillations in zones, which were not recommended for operation,
were found. The obtained results are consistent with the full-scale test results carried out
for hydropower units of the Sayano-Shushenskaya hydropower plant in 1988.
\\
\end{abstract}

\section{Introduction}

Nowadays one of the most important source of electricity is the electric energy produced by hydroelectric facilities.
According to \cite{Wagner2011}, Norway gets 99\% of its electric power from water, Brazil 84\%, Austria 59\%, Canada 58\% and Russia 18\%.
Failures of hydropower plants cause loss of lives, major damage in the surrounding area, and serious economic consequences.
In recent years,  accidents at the hydropower plants have become frequent (see, e.g., Bieudron Hydroelectric Power Station (Switzerland, 2000),
Taum Sauk Hydroelectric Power Station (Missouri, USA, 2005), Sayano–Shushens\-kaya Dam (Russia, 2009), Itaipu Dam (Brazil, 2009), Srisailam Dam (India, 2013), Dhauliganga hydro electric station (India, 2013)). In order to prevent such accidents, it is necessary to investigate their causes.

\begin{figure}[!ht]
	\centering 		
		\includegraphics[width=1\linewidth]{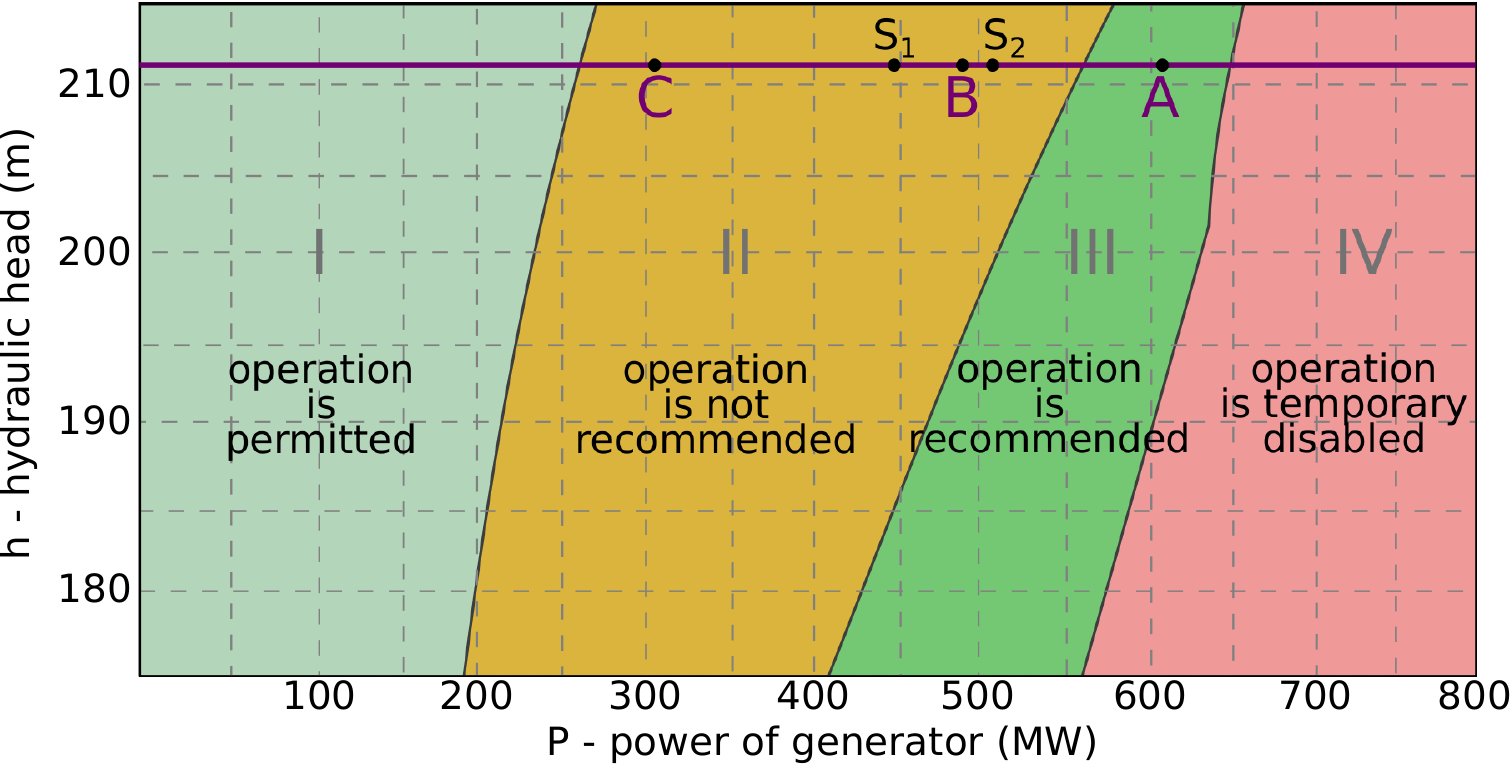}
		\caption{Operational zones of hydropower unit of the Sayano-Shushenskaya hydropower plant}
	\label{zones_1}
\end{figure}

This work is motivated by the accident happened on the Sayano-Shushenskaya hydropower plant in 2009 year.
According to the act of special commission of \emph{the Russian Federal Environmental, Industrial and Nuclear Supervision Service}, 
immediately before the accident the power of the second hydropower unit was 475~MW at a head of 212 meters \cite{Act}, i.e., it worked in the not recommended zone II (Fig.~\ref{zones_1}).
Zone II of hydropower unit work is characterized by strong hydraulic turbine blows in flowing part and vibrations.
For the normal operation it is recommended power range, corresponding to the zone III, in which the efficiency of
turbines has a maximum value. Also operation is allowed in the zone I, in which the dynamics is allowed, 
but the level of efficiency of the turbines are low. Operation in zone IV is not allowed. 
These work zones of hydropower unit of the Sayano-Shushenskaya hydropower plant were obtained by the full-scale test of hydropower units 
in the late 80s of 20th century and were published in the
technical report ``{\it Full-scale testing of turbines of Sayano-Shushenskaya hydropower plant with standard runner}'' No. 1008 \cite{Gluhih}.

It is known \cite{Act}, that the accident was caused by vibrations in the hydropower unit. Various possible reasons of occurrence of such vibrations were discussed in \cite{Act,Kurzin}. However main attention in these works was given to the dynamical processes in the turbines 
and the full mathematical model turbine--speed governor--synchronous generator was not considered.

In this paper an adequate mathematical model turbine--speed governor--synchronous generator is suggested. 
This model joins together differential equations of turbine (Pervozvanskiy's model) \cite{Pervozvanskii}, differential equations of synchronous generator (Park-Gorev model) \cite{Adkins,boldea}, and differential equations of speed governor \cite{Merkurev}. 
Parameters of the model are taken from \cite{Karapetyan}. 
As a result of using numerical methods the periodic solutions of differential equations of the model were found during the analysis of transient processes.
These solutions correspond to hydropower unit vibrations observed for some operation modes at the Sayano-Shushenskaya hydropower plant.

Remark, that we considered also some other more simple models of hydropower unit 
(based on more simple models of turbines (see \cite{Merkurev}), synchronous generators (the Tricomi equation \cite{1931-Tricomi,leonov2009},
the equations of synchronous generator with parallel connection in feed system \cite{leonov2009}).
However, for these models the above effects have not been found.

\section{Components of hydropower unit}
\begin{figure}[!ht]
	\centering 		
		\includegraphics[width=0.7\linewidth]{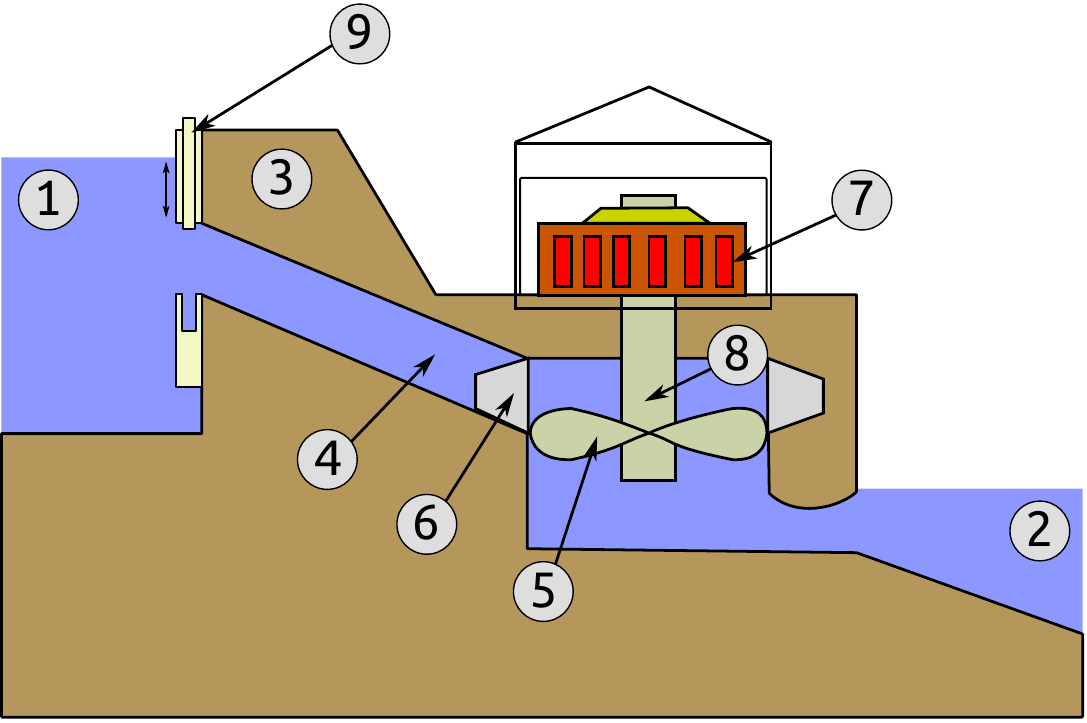}
		\caption{Scheme of the Sayano-Shushenskaya hydropower plant: 1 -- upper reservoir, 2 -- lower reservoir, 3 -- dam, 4 -- penstock, 5 --  blade wheel, 6 -- guide vanes, 7 -- rotor of synchronous generator, 8 -- shaft, 9 -- control gate}
		\label{construction}
\end{figure}

Following \cite{LeonovKS-2015}, let us consider the main constructive elements of hydropower unit. 
The hydropower unit in Fig.~\ref{construction} consists 
of synchronous generator and radial-axial hydraulic turbine. 
In this work a nonregulated synchronous generator is considered. 
The rotor of generator and the runner are connected together by a rigid shaft.

The dam creates a difference in water level between the upper reservoir and lower reservoir. The flow of water is delivered from  upper reservoir to turbine by penstock through the spiral casing. 
Water jets impact on the blades of the turbine producing torque applied to the rotating shaft.
Since the turbine shaft is rigidly connected with the generator rotor,
the rotor starts to rotate and to produce electricity, which is transferred to the grid.
The water flow is controlled by means of guide vanes.

For the safety of the power network, the frequency should remain almost constant. 
This is reached by keeping the same speed of the synchronous generator. 
The rotational speed is controlled by the speed governor.

\begin{figure}[!ht]
	\centering 		
		\includegraphics[width=1\linewidth]{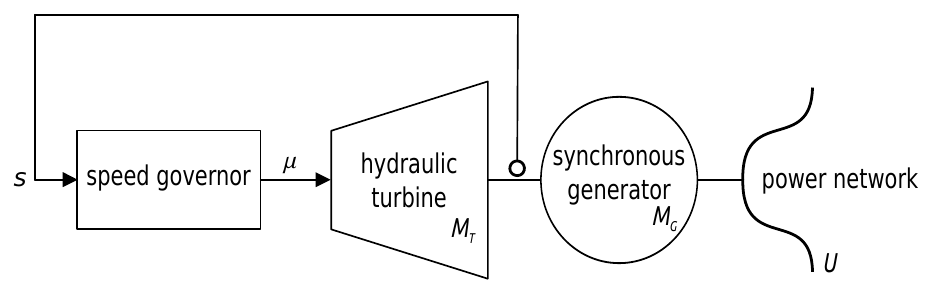}
		\caption{The main structural elements of hydropower unit
		}
		\label{scheme}
\end{figure}

The main structural elements of hydropower unit are presented in Fig.~\ref{scheme}.
Introduce the following notations:
	$s$  is a signal of angular (rotational) speed deviation,
	$\mu$ is a position of guide vanes,
	$M_T$ is a turbine torque,
	$M_G$ is a generator torque,
    $U$  is  a voltage required by the power network.
The power network represents a set of energy suppliers and consumers.

Thus, the hydropower unit represents a system consisting of hydraulic turbine, which produces a mechanical torque, generator,
which converts mechanical energy into electrical energy, and automatic speed governor, 
which regulates the rotation speed of hydraulic turbine.

In order to develop a mathematical model of hydropower unit it is necessary to describe each structural element of
hydropower unit, presented in Fig.~\ref{scheme}. 
Let us consider each elements in details.

\section{Mathematical model of synchronous generator} \label{sec_generator}

\begin{figure}[!ht]
	\centering 		
		\includegraphics[width=0.40\linewidth]{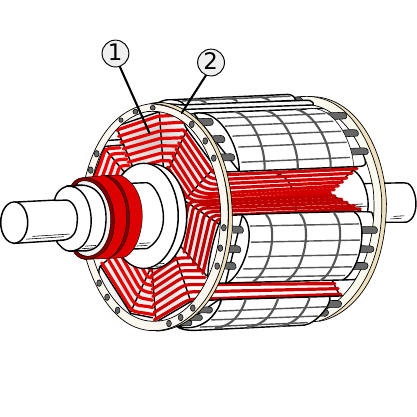} \hfill
		\includegraphics[width=0.35\linewidth]{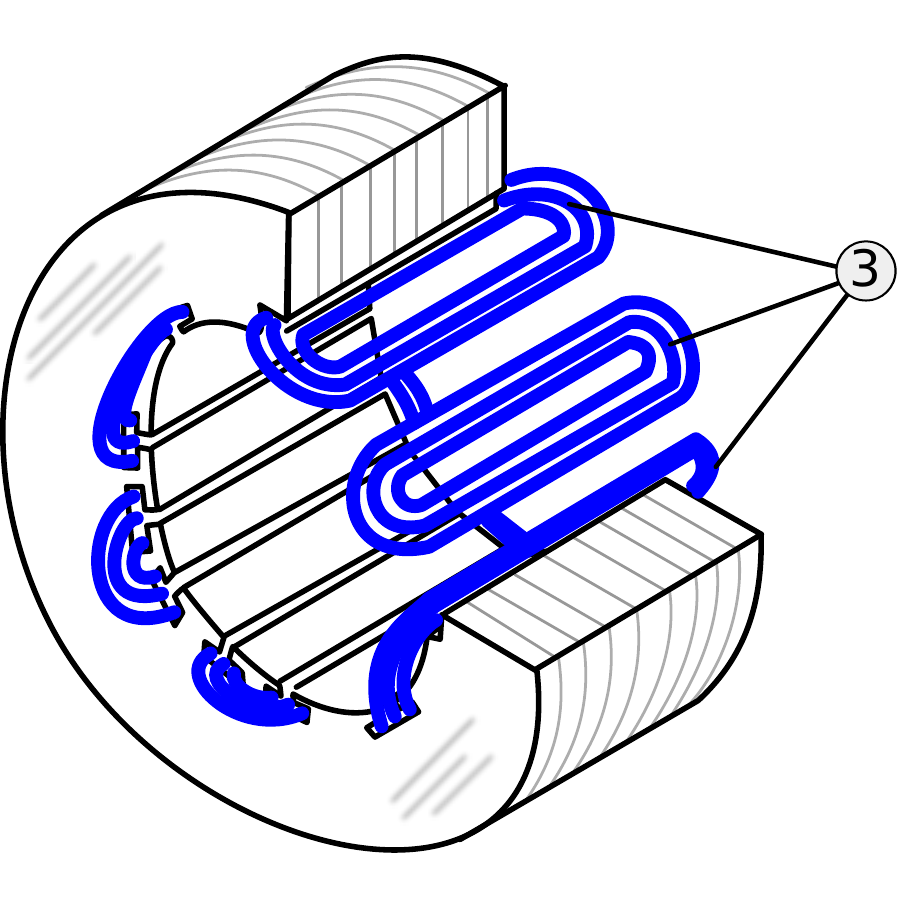} \\
		a \hspace*{0.49\linewidth} b
		\caption{Synchronous generator. (a) -- rotor: 1 -- field winding, 2 -- damper winding; 
                                        (b) -- stator: 3 -- stator winding}
		\label{pic_rotor_stator}
\end{figure}

The main constructive elements of synchronous generator are stationary stator (Fig.~\ref{pic_rotor_stator}, a) and rotating rotor (Fig.~\ref{pic_rotor_stator}, b). The windings are placed in the stator and rotor slots. Stator windings are arranged in such a way that
in the case of alternate current they generates a rotating magnetic field.

Using  Park's transformation, the three-phase windings of  stator can be substituted by two equivalent short-circuited windings, and the
rotor can be described as two equivalent short-circuited damper windings and one field winding. 
Let us choose the rotating coordinate system ($ d $,~$ q $), connected with the rotor.
\begin{figure}[!ht]
	\centering 		
		\includegraphics[width=0.6\linewidth]{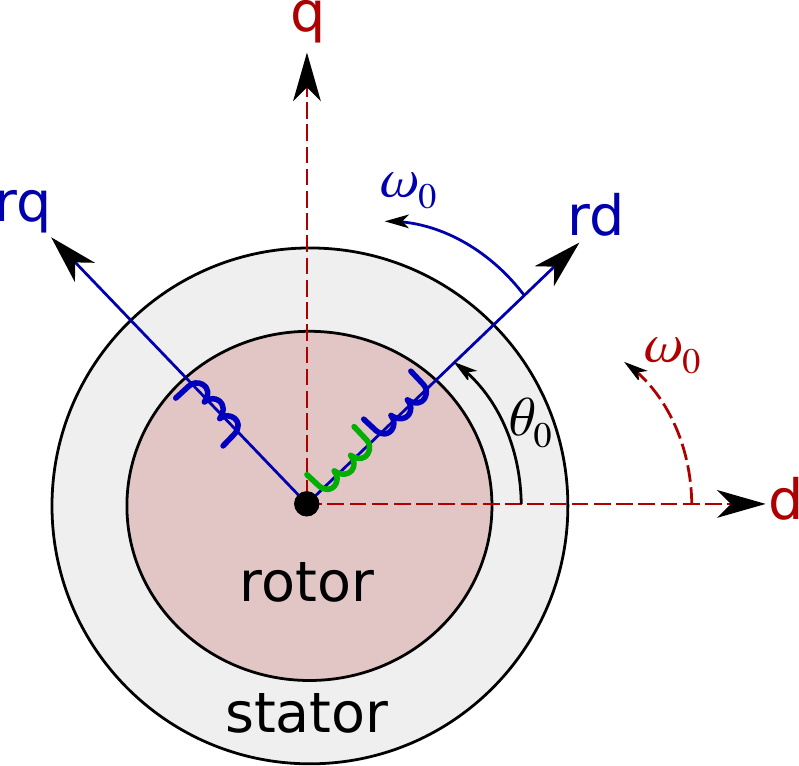}
		\caption{The rotating coordinate system ($ d $,~$ q $): axes ($rd$, $rq$) are rigidly connected with rotating rotor;
		axes ($d$, $q$) are rigidly connected with rotating magnetic field of stator}
		\label{scheme-coordinates}
\end{figure}

Introduce the following notations:
$$
\theta(t) = \theta_0 + \theta_\Delta(t),
$$
$$
\omega(t) = \omega_0 + \omega_\Delta(t),
$$
$$
s(t) = \frac{\omega_\Delta(t)}{\omega_0},
$$
where  $\theta$ is an angle between rotor axes and rotating magnetic field, generated by the stator windings [rad],
$\theta_0$ is an operating electrical angle [rad] ($0 \leq \theta_0 \leq \pi/2$),
$\theta_\Delta$ is a deviation of operating angle,
$\omega$ is an angular speed of rotor (shaft) [rad/s],
$\omega_0$ is a rated angular speed (synchronous speed) [rad/s],
$\omega_\Delta$ is a deviation of rated angular speed,
$s$ is a relative deviation of rated angular speed [p.u.].
The rated angular speed is determined by the formula
$$
\omega_0 = \frac{2\pi f }{p},
$$
where $f = 50$ is the network frequency [$\rm{Hz} = s^{-1}$],
$p$ is the number of pairs of poles.
Since hydrogenerators of the Sayano-Shushenskaya hydropower plant have $21$ pairs of poles,
then $\omega_0 = 2 \pi 50 / 21 = 2\pi \cdot 2.38$ [rad/s] or $\omega_0 = 2 \pi 3000 / 21 = 2\pi \cdot 142.85$ [rad/min].
Synchronous generators are usually operated with the nominal power factor $ \cos \theta_0 = 0,8 ... 0,9 $ \cite{Karapetyan}.
For generators of the Sayano-Shushenskaya hydropower plant $\theta_0 = \arccos (0,9)$.

Thus, a synchronous generator can be described by the Park-Gorev equations in per unit values \cite[pp.132-149]{Adkins}, \cite[pp.5.15-
5.17]{boldea}:\\
\begin{equation} \label{equation_PG}
\begin{split} &
	\dot\Psi_d =
				- \omega_0 (1 + s) \Psi_q - \omega_0\,r\,i_d - \omega_0\,u_d,
\\ &
	\dot\Psi_q =
				\omega_0 (1 + s) \Psi_d - \omega_0\,r\,i_q - \omega_0\,u_q,
\\ &
	\dot\Psi_r       = \frac{1}{T_r} (E_r - E_q),
\\ &
	\dot\Psi_{rd} = -\frac{1}{T_{rd}} E_{rq},
\\ &
	\dot\Psi_{rq} = \frac{1}{T_{rq}} E_{rd},
\end{split}
\end{equation}
where flux linkages are determined as follows
\begin{equation} \label{equation_flux}
\begin{split} &
	\Psi_d    = x_d\,i_d + E_q + E_{rq},
\\ &
	\Psi_q    = x_q\,i_q - E_{rd},
\\ &
	\Psi_r    = \frac{x_{ad}^2}{x_r}i_d + E_q + \frac{x_{ad}}{x_r}E_{rq},
\\ &
	\Psi_{rd} = \frac{x_{ad}^2}{x_{rd}} i_d + E_{rq} + \frac{x_{ad}}{x_{rd}} \left( E_q + E_{rq} \right),
\\ &
	\Psi_{rq} = \frac{x_{aq}^2}{x_{rq}} i_q - E_{rd} - \frac{x_{aq}}{x_{rq}} E_{rd}.
\end{split}
\end{equation}

Here the following variables and coefficients are relative 
to the corresponding base values (voltage, current, flux linkage, impedance, inductance, power):
$r$ is a stator resistance,
$i_d$, $i_q$ are currents in stator windings,
$u_d$, $u_q$ are  stator voltages,
$\Psi_d$, $\Psi_q$, $\Psi_r$, $\Psi_{rd}$, $\Psi_{rq}$ are flux linkages,
$E_r$ is a field voltage,
$E_q, E_{rd}, E_{rq}$ are electromotive forces, induced in the stator by
the magnetic field of rotor winding currents for synchronous rotor speed,
$x_d$, $x_q$ are synchronous inductances (reactances) along the axes $d$ and $q$,
$x_r$, $x_{rd}$, $x_{rq}$ are impedances of field winding, damper
windings along the axes $d$ and $q$,
$x_{ad}$, $x_{aq}$ are impedances of stator winding along the axes $d$ and $q$.
The following coefficients are time constants:
$T_r$ is a field-winding time constant with open stator and damper windings [s],
$T_{rd}, T_{rq}$ are damper winding time constants with open stator and field windings [s].

The stator voltages $u_d$ and $u_q$ along the d- and q- axes in per unit
values are determined according to the laws
$$u_d = - U \sin(\theta_0+\theta_\Delta), \qquad  u_q = U \cos(\theta_0+\theta_\Delta),$$
where $U = \frac{\widehat{U}}{U_\text{b}}$ is a voltage relative to the
base voltage [p.u.], $\widehat{U}$ is a voltage in power network [V],
$U_\text{b}$ is the base voltage [V].

The motion of synchronous generator rotor about shaft is described by the torque equation in physical values \cite[pp.132-149]{Adkins},
\cite[pp.5.15-5.17]{boldea}:
\begin{equation} \label{eq_motion}
\begin{split} &
	\dot\theta_\Delta = \omega,
\\ &
	J\,\dot{\omega} =  M_T - \left( \widehat{\Psi}_d \, \widehat{i}_q -
\widehat{\Psi}_q \, \widehat{i}_d \right).
\end{split}
\end{equation}
Here
$J$ is a moment of inertia of rotor 
(is a moment of inertia of hydropower unit) [$\text{kg}\cdot \text{m}^2$], $\widehat{\Psi}_d$,
$\widehat{\Psi}_q$ are flux linkages in physical unit values [$\rm{Wb} =
\frac{\text{kg}\cdot \text{m}^2}{\text{s}^2\text{A}}$],
$\widehat{i}_d$, $\widehat{i}_q$ are currents in stator windings in
physical unit values [$\text{A}$],
$M_T$ is the rotation torque [$\frac{\text{kg}\cdot
\text{m}^2}{\text{s}^2}$].
Since the shaft of turbine is rigidly connected with the rotor of
generator, a moment of inertia of rotor coincides with a moment of
inertia of hydropower unit and an angular rotor speed coincides with an
angular turbine speed. In our case the rotation torque $M_T$ is  a
turbine torque, which is created by the pressure of water on the runner.

Rewrite system of equations~(\ref{eq_motion})  in terms of base variables.
For this purpose the second equation is divided twice by the base voltage:
$$
U_\text{b} = \omega_0 \Psi_\text{b} = Z_\text{b} I_\text{b},
$$
where $\Psi_\text{b}$ is the base flux linkage [$\rm{Wb} =
\frac{\text{kg}\cdot \text{m}^2}{\text{s}^2\text{A}}$], $Z_\text{b}$ is
the base impedance (valid also for resistances and reactances)[Ohms],
$I_\text{b}$ is the base current [A].
Then
$$
\frac{J}{\Psi_\text{b} Z_\text{b} I_\text{b}}\,\left(
\dot{\frac{\omega}{\omega_0}} \right) =  \frac{M_T}{\omega_0
\Psi_\text{b} Z_\text{b} I_\text{b}} - \frac{1}{\omega_0 Z_\text{b}}
\left( \frac{\widehat{\Psi}_d}{\Psi_\text{b}} \,
\frac{\widehat{i}_q}{I_\text{b}} - \frac{\widehat{\Psi}_q}{\Psi_\text{b}}
\, \frac{\widehat{i}_d}{I_\text{b}} \right),
$$
$$
\frac{J\omega_0^2}{U_\text{b} I_\text{b}}\,\dot{s} =
\frac{M_T}{\Psi_\text{b} I_\text{b}} -\left( \Psi_d \, i_q - \Psi_q \,
i_d \right).
$$

Thus, as a result of given transformation one obtains system of equations~(\ref{eq_motion}) in per unit values:
\begin{equation} \label{eq_moment}
\begin{split} &
	\dot\theta_\Delta = \omega_0 s,
\\ &
	\dot{s} = \frac{1}{T_J} \left(\frac{\widehat{M_T}}{\Psi_\text{b}\,I_\text{b}} - \left( \Psi_d \, i_q - \Psi_q \, i_d \right)\right),
\end{split}
\end{equation}
where $T_J = J \frac{\omega_0^2}{S_\text{b}}$ is inertial constant of hydropower unit~[s],
$S_\text{b} = U_\text{b} I_\text{b} $ is the base power [$\text{W} = \frac{\text{kg}\cdot\text{m}^2}{\text{s}^3}$].

\section{Mathematical model of hydraulic turbine}

Hydraulic turbines can be classified by their type of construction, the most important ones being the Francis, Pelton and Kaplan or
Propeller turbines \cite{Wagner2011}. They are distinguished by construction of runner and control methods of the rotation speed of the turbine. However, regardless of the turbine structure, the rotation equation of turbine shaft in earth-fixed coordinate system has the following form \cite{kundur1994power,2001-Springer-Leonhard}
\begin{equation*}
\begin{split} &
    	J\dot \omega = M_T - M_G,
	\end{split}
\end{equation*}
where
$M_G$ is a resistance torque, which is created by the generator.

The rotational torque is defined by the relation
\begin{equation*}
	\begin{split} &
		M_T = \widetilde{k} \frac{P_T(t)}{\omega(t)},
	\end{split}
\end{equation*}
where $P_T$ is the turbine mechanical power [$\text{W} = \frac{\text{kg}\cdot\text{m}^2}{\text{s}^3}$],
$\widetilde{k}$ is a constant depending on the construction of turbine [p.u.].

By virtue of the laws of hydraulics, the turbine power $ P_T $ is
calculated by the formula:
\begin{equation*}
	\begin{split} &
		P_T (t) = h(t) Q(t),
	\end{split}
\end{equation*}
where $h$ is a pressure drop
[$\text{Pa}=\frac{\text{kg}}{\text{m}\text{s}^2}$], $Q$ is the water flow
through the turbine [$\frac{\text{m}^3}{s}$].

The water flow through the turbine is related with $h$ and the opening
(closing) guide vanes in the following way:
$$
Q(t) = {\rm C} \mu(t) \sqrt{h(t)},
$$
$$
\mu (t) = \mu_0 + \mu_\Delta (t),
$$
where
$\mu$ is a position of guide vanes [p.u.],
$\mu_0$ is a given position of guide vanes [p.u.],
$\mu_\Delta$ is the deviation of given position of guide vanes [p.u.],
${\rm C} = S/\sqrt{\rho}$  is a constant depending on the construction of
 penstock [$\frac{\text{m}^3\sqrt{\text{m}}}{\sqrt{\text{kg}}}$],
$S$ is a sectional area of water conduit [$\text{m}^2$],
$\rho$ is the density of water [$\frac{\text{kg}}{\text{m}^3}$].

The differential equation, describing the change of water flow
through the turbine, can be written according to the Newton equation
for the liquid column enclosed in the penstock. Then, taking into account that
$S$ is a sectional area of penstock [$\text{m}^2$],
$l$  is a  length of penstock [$\text{m}$],
$\rho$ is  a density of water [$\frac{\text{kg}}{\text{m}^3}$],
$p_u$ is a constant pressure on the upper end of penstock
[$\text{Pa}$],
$p_l$ is  a constant pressure on the lower end of penstock (after
turbine) [$\text{Pa}$],
 one obtains
\begin{equation*}
	\begin{split} &
		\dot Q= \frac{S}{l\rho}(p_u-p_l-h(t)).
	\end{split}
\end{equation*}

Thus, the dynamics of hydraulic turbine is described by the system of
differential equations
\begin{equation*}
	\begin{split} &
		J \,\dot \omega =
				\frac{\widetilde{k}}{C^2(\mu_0+\mu_\Delta)^2} \frac{Q^3}{\omega} - M_G,
	\\ &
		\dot Q = \frac{S}{l\rho} \left(p_u - p_l - \frac{Q^2}{C^2(\mu_0+\mu_\Delta)^2} \right),
	\end{split}
\end{equation*}
where $M_G$ is a generator torque in the physical values.

In order to pass to the rotating coordinates ($d,q$), 
it is necessary to make the following transformation
\begin{equation} \label{eq_preob}
 \omega_\Delta(t) = \omega(t) - \omega_0.
\end{equation}

Since for the description of dynamics of synchronous generator 
there are used equations in per unit values, 
we write the motion equation of turbine also in per unit values:
\begin{equation*}
	\begin{split} &
		\dot s = \frac{1}{T_J} \left(
				\frac{k}{C^2(\mu_0+\mu_\Delta)^2} \frac{Q^3}{\omega_0^2 (1 + s)} - \widetilde{M_G} \right),
	\\ &
	\end{split}
\end{equation*}
where
$$
k =  \frac{\widetilde{k}}{\Psi_\text{b}\,I_\text{b}}
$$
and in  Section~\ref{sec_generator} there was defined
the generator torque in per unit values
$$\widetilde{M}_G =  \Psi_d \, i_q - \Psi_q \, i_d .$$

Note that the obtained system contains one control signal $\mu_\Delta$, which corresponds to  automatic speed governor of  turbine.

\section{Mathematical model of speed governor} \label{section-turbine-law}

\begin{figure}[ht!]
\includegraphics[width=1\linewidth]{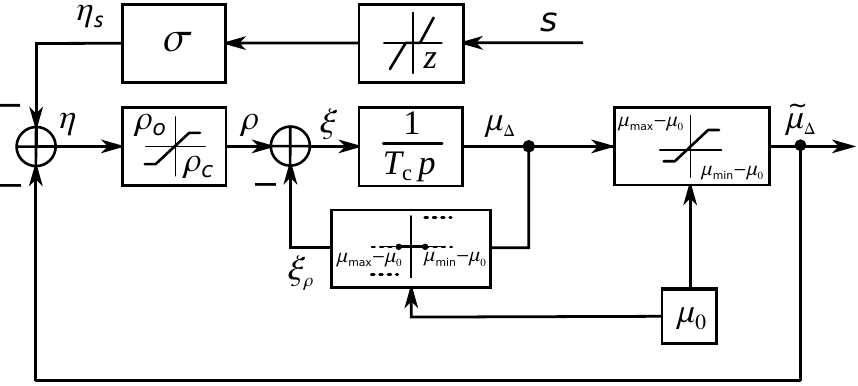}
\caption{Scheme of automatic speed governor of hydraulic turbine}
\label{turbina}	
\end{figure}

Let us consider the simplified scheme of automatic speed governor of turbine \cite[pp.153-157]{Merkurev}, presented in Fig.~\ref{turbina}.
In this work parameters of hydropower unit are chosen so that the speed governor does not get to saturation. 
To simulate the case of saturation 
it is necessary to consider more complex models of speed governors
(see, e.g., \cite{bondareva2004,nicolet2007,sattoufsimulation2014,TarbouriechGarcia_book11,TarbouriechTurner09}).

The governor has the deadband $z$, which is specified by the technical conditions. 
The deadband of the hydraulic turbine is $30 \, \rm{mHz}$ \cite{Merkurev}. 
The deadband in the per unit values is $z = 0,002$.

The input signal is a relative deviation of rated angular speed $s$.
After the input signal passes through the deadband, one obtains the signal
$ \eta_s $, which corresponds to  a signal of measuring device:
$$
		\eta_s = \sigma\chi_s(s) = \left\{
			\begin{array}{lc}
				\sigma (s - z/2), & s \geq z/2,
			\cr
				\sigma (s + z/2), & s \leq -z/2,
			\cr
				0,  & \lvert s \rvert < z/2,
			\end{array}
		\right.
$$
where $\sigma$ is a transmission coefficient of open-cycle control system.

The control signal $ \eta $ is formed by the formula
$$
\eta = - \eta_s - \widetilde{\mu}_\Delta,
$$
where $\widetilde{\mu}_\Delta$  is  a signal of rigid negative feedback.

Then the control signal is cut off because of restriction on the velocity
of change of guide vanes position:
$$
\rho = \chi_\rho (\eta) = \left\{
			\begin{array}{lc}
				\eta, & \rho_o \leq \eta \leq \rho_c,
			\cr
				\rho_o, & \eta < \rho_o,
			\cr
				\rho_c, & \rho_c < \eta,
			\end{array}
		\right.
$$
where $\rho_o, \rho_c$ are the maximum velocities of opening and closing the vanes.

A servomotor is presented by the integrator with time constant  $T_c$ [s]. The output value of servomotor is the relative displacement of guide vanes $\widetilde{\mu}_\Delta$.
The stroke of servomotor has limit stops $\mu_{min}, \mu_{max}$, which correspond to the minimum and maximum power of the turbine \cite{Merkurev}. If guide vanes reach limit stops, at which further displacement of vanes in the same direction is not possible, then displacement of vanes is stopped. 
In some research works it is recommended to model such saturation through direct feedback \cite{GaleaniOnoriTeel_SCL08,TarbouriechGarcia_book11,TarbouriechTurner09,ZaccarianTeel_11book}.

Consequently, the control signal takes the following form
$$
	\xi = \rho - \xi_\rho(\mu_\Delta, \mu_0),
$$
where
$$
\xi_\rho(\mu_\Delta, \mu_0) = \left\{
			\begin{array}{lcc}
				\rho, & & \mu_\Delta + \mu_0 < \mu_{min} \text{ and } \rho < 0,
			\cr
				\rho, & & \mu_\Delta + \mu_0 > \mu_{max} \text{ and } \rho > 0,
			\cr
				0, & & \text{otherwise}.
			\end{array}
		\right.
$$
The equation of servomotor motion is as follows
\begin{equation}
	\dot{\widetilde{\mu}}_\Delta= \frac{\rho - \xi_\rho(\mu_\Delta, \mu_0)}{T_c},
	\label{eq_servo}
\end{equation}
where\\
$
\widetilde{\mu}_\Delta = \chi_\mu(\mu_\Delta, \mu_0) =
$
$$
 \text{~~~~~~~} =\left\{
			\begin{array}{lc}
				\mu_\Delta, & \mu_{min} -\mu_0 \leq \mu_\Delta \leq \mu_{max} -\mu_0,
			\cr
				\mu_{min} - \mu_0, & \mu_0 + \mu_\Delta < \mu_{min},
			\cr
				\mu_{max} - \mu_0, & \mu_0 + \mu_\Delta < \mu.
			\end{array}
		\right.
$$

In other words, the signal $\widetilde{\mu}_\Delta$ is the signal $\mu_\Delta$, passed through the saturation. Consequently, equation (\ref{eq_servo}) can be rewritten in the form
$$
	\dot\mu_\Delta = \frac{\rho - \xi_\rho(\mu_\Delta, \mu_0)}{T_c},
$$
$$
\mu_{min} -\mu_0 \leq \mu_\Delta \leq \mu_{max} -\mu_0.
$$

Thus, the automatic speed governor of the turbine can be
described by the following differential equation
$$
	\begin{array}{l} \displaystyle
			\dot\mu_\Delta = \frac{
				\chi_\rho\left(- \sigma \chi_s(s) - \chi_\mu(\mu_\Delta, \mu_0) - \xi_\rho(\mu_\Delta, \mu_0) \right)
			}{T_c}.
		\end{array}
$$

\section{Mathematical model of hydropower unit. Operating modes}

Using equations of each structural element of hydropower unit,
one writes the equations of hydropower unit with automatic speed governor in per unit values:
\begin{equation} \label{model}
\begin{split} &
		\dot\theta_\Delta = \omega_0\,s,
\\ &
	\dot s =
				\frac{1}{T_J}\Big(\frac{k}{C^2(\mu_0+\mu_\Delta)^2} \frac{Q^3}{\omega_0^2 (1 + s)} -  \Psi_d i_q + \Psi_q i_d \Big),
\\ &
		\dot Q = \frac{S}{l\rho} \left(p_u - p_l - \frac{Q^2}{C^2(\mu_0+\mu_\Delta)^2} \right),
\\ &
	\dot\Psi_d =
				- \omega_0 (1 + s) \Psi_q - \omega_0\,r\,i_d + \omega_0\,U \sin (\theta_0 + \theta_\Delta),
\\ &
	\dot\Psi_q =
				\omega_0 (1 + s) \Psi_d - \omega_0\,r\,i_q - \omega_0\,U \cos (\theta_0 + \theta_\Delta),
\\ &
	T_r\,\dot\Psi_r       = E_r - E_q,
\\ &
	T_{rd}\,\dot\Psi_{rd} = - E_{rq},
\\ &
	T_{rq}\,\dot\Psi_{rq} = E_{rd},
\\ &
	\dot\mu_\Delta = \frac{
				\chi_\rho\left(- \sigma \chi_s(s) - \chi_\mu(\mu_\Delta, \mu_0) - \xi_\rho(\mu_\Delta, \mu_0) \right)
			}{T_c},
\\&
	\Psi_d    = x_d\,i_d + E_q + E_{rq},
\\ &
	\Psi_q    = x_q\,i_q - E_{rd},
\\ &
	\Psi_r    = \frac{x_{ad}^2}{x_r}i_d + E_q + \frac{x_{ad}}{x_r}E_{rq},
\\ &
	\Psi_{rd} = \frac{x_{ad}^2}{x_{rd}} i_d + E_{rq} + \frac{x_{ad}}{x_{rd}} \left( E_q + E_{rq} \right),
\\ &
	\Psi_{rq} = \frac{x_{aq}^2}{x_{rq}} i_q - E_{rd} - \frac{x_{aq}}{x_{rq}} E_{rd}.
\end{split}
\end{equation}

In the case of the steady-state stability of hydropower unit the turbine rotates with a constant angular speed ($s=0$) and the generator produces the constant power. Such an operation mode of hydropower unit corresponds to a steady-state (operating) mode of power network.
Note that for computations
the hydropower unit are often presented as a source of current or electromotive force,
the power of which corresponds to a power of generator.

Dynamical stability of hydropower units is considered in terms of maintaining a given mode. This means that if sudden, significant changes of  network mode arises, then after the transient processes the output power of hydropower unit must correspond to the required power. For example, after short circuits in one or more power lines, blackouts, changes of the external load, etc.
Note that the hydropower unit in the considered processes
can not be represented as a source of current or electromotive force
since dynamical processes in this case have a significant effect both on a hydropower unit, and on the network mode.

An operating mode of hydropower unit corresponds to the asymptotically stable equilibrium point of system (\ref{model}).
The equilibrium points are the following
\begin{equation*}
\begin{split} &
	s^{\text{st}}=0  \quad (\text{i.e.} \quad \omega_\Delta^{\text{st}} = 0), \qquad	\mu_\Delta^{\text{st}}=0, \qquad
\\ &
	\theta_\Delta^{\text{st}}= \rm{const},  \qquad	Q^{\text{st}}=C\mu_0\sqrt{p_u-p_l},
\\&
E^{\text{st}}_q=E_r, \qquad E^{\text{st}}_{rq}=0, \qquad 	E^{\text{st}}_{rd}=0,
\\ &
	i^{\text{st}}_d= - \frac{x_q}{r^2+x_dx_q} (-\frac{r}{x_q}U \sin \theta - U \cos \theta + E_r),
\\&
	i^{\text{st}}_q= - \frac{r}{r^2+x_dx_q} (-\frac{x_d}{r}U \sin \theta - U \cos \theta - E_r),
\\ &
	\Psi^{\text{st}}_d= x_di_d+E_r,  \qquad \Psi^{\text{st}}_q=x_qi_q, \qquad \Psi^{\text{st}}_r=\frac{x_{ad}^2}{x_r}i_d+E_r,
\\ &
		\Psi^{\text{st}}_{rd}=\frac{x_{ad}^2}{x_{rd}}i_d+\frac{x_{ad}}{x_{rd}}E_r, \qquad \Psi^{\text{st}}_{rq}=\frac{x_{aq}^2}{x_{rd}}i_q.
\end{split}
\end{equation*}
The position of guide vanes $\mu_0$ is defined from the following equation, which will be called the balance equation between the turbine and the
generator torques:
\begin{equation} \label{balance}
	M_T(\mu_0, \omega_0) = M_G(U,\theta_0),
\end{equation}
where
$$
	M_T(\mu_0, \omega_0) = \frac{kC(p_u - p_l)^\frac{3}{2}}{\omega_0^2} \mu_0,
$$
$$
	\begin{array}{l} \displaystyle
		M_G (U,\theta)=
			\frac{r(x_d-x_q)}{r^2+x_dx_q)^2} (
				- U^2x_d\sin^2\theta + U^2x_q\cos^2\theta
				-
	\cr\displaystyle
			- \frac{r^2-x_dx_q}{r}U^2 \sin\theta\cos\theta
			- E_r\frac{x_dx_q-r^2}{r}U\sin\theta
			+
	\cr\displaystyle
			+ 2x_qE_rU\cos\theta) - \frac{E_r}{r^2+x_dx_q}
				(-x_dU\sin\theta + rU\cos\theta)
			+
	\cr\displaystyle
			+ \frac{rx_qE_r(x_d-x_q)}{(r^2+x_dx_q)^2}
			- \frac{rE_r^2}{(r^2+x_qx_d)}.
	\end{array}
$$

\begin{figure}[!ht]
	\begin{center}
		\includegraphics[width=1\linewidth]{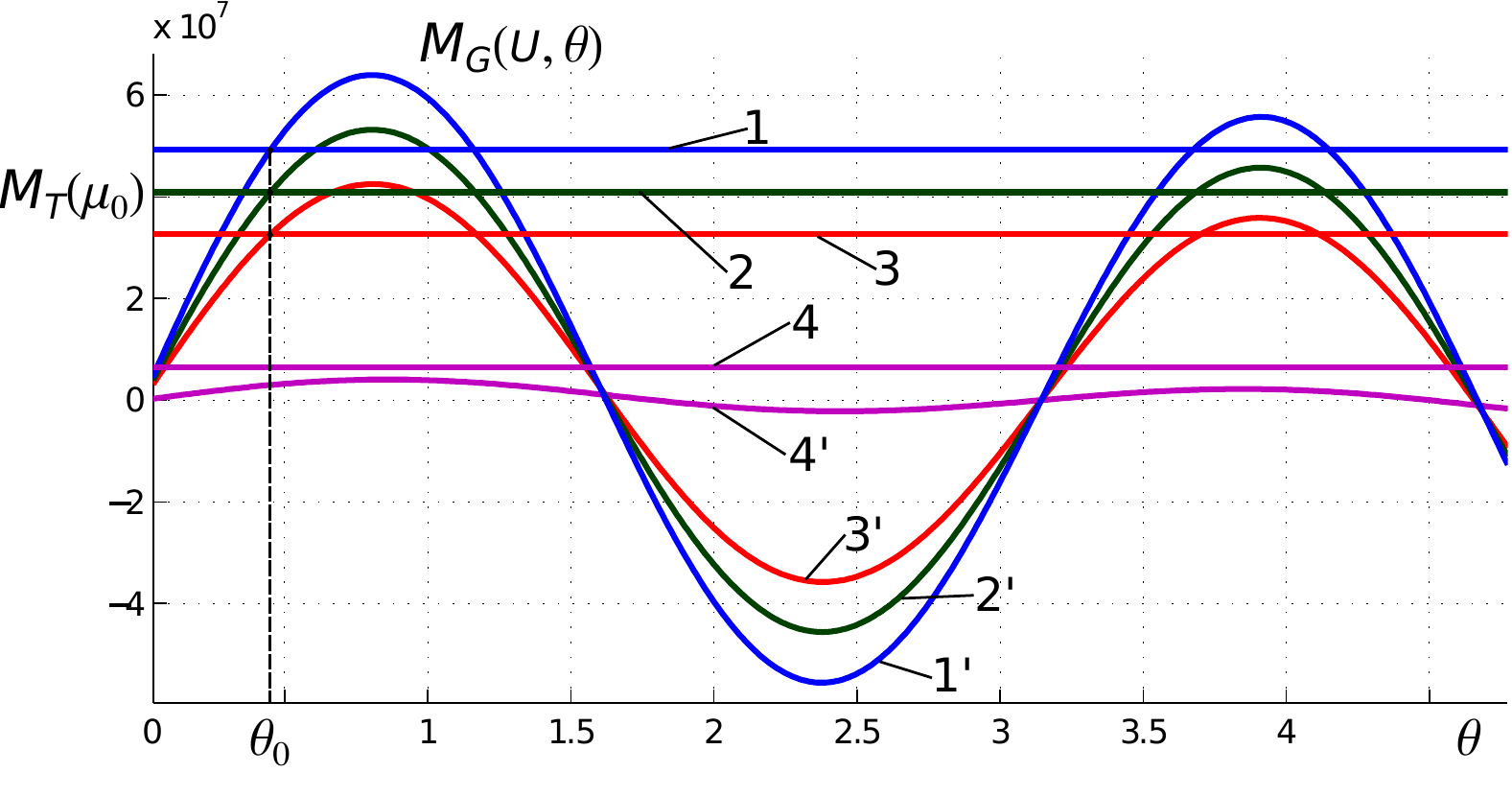}
	\caption{Graphical solution of the balance equation depending on $U$.
     Lines 1, 1' for $\gamma = 1$, lines 2, 2' for $\gamma = 0.9$,
	 lines 3, 3' for $\gamma = 0.8$, lines 4, 4' for $\gamma = 0.1$
}
	\label{pic_statics}
	\end{center}
\end{figure}

Graphical solutions of the balance equation depending on voltage change in the power network are presented in Fig. \ref{pic_statics}.
The parameter corresponding to the voltage is represented as $U = \gamma U_{nom}$, where $\gamma  > 0$.

\begin{figure}[!ht]
\centering
	\includegraphics[width=1\linewidth]{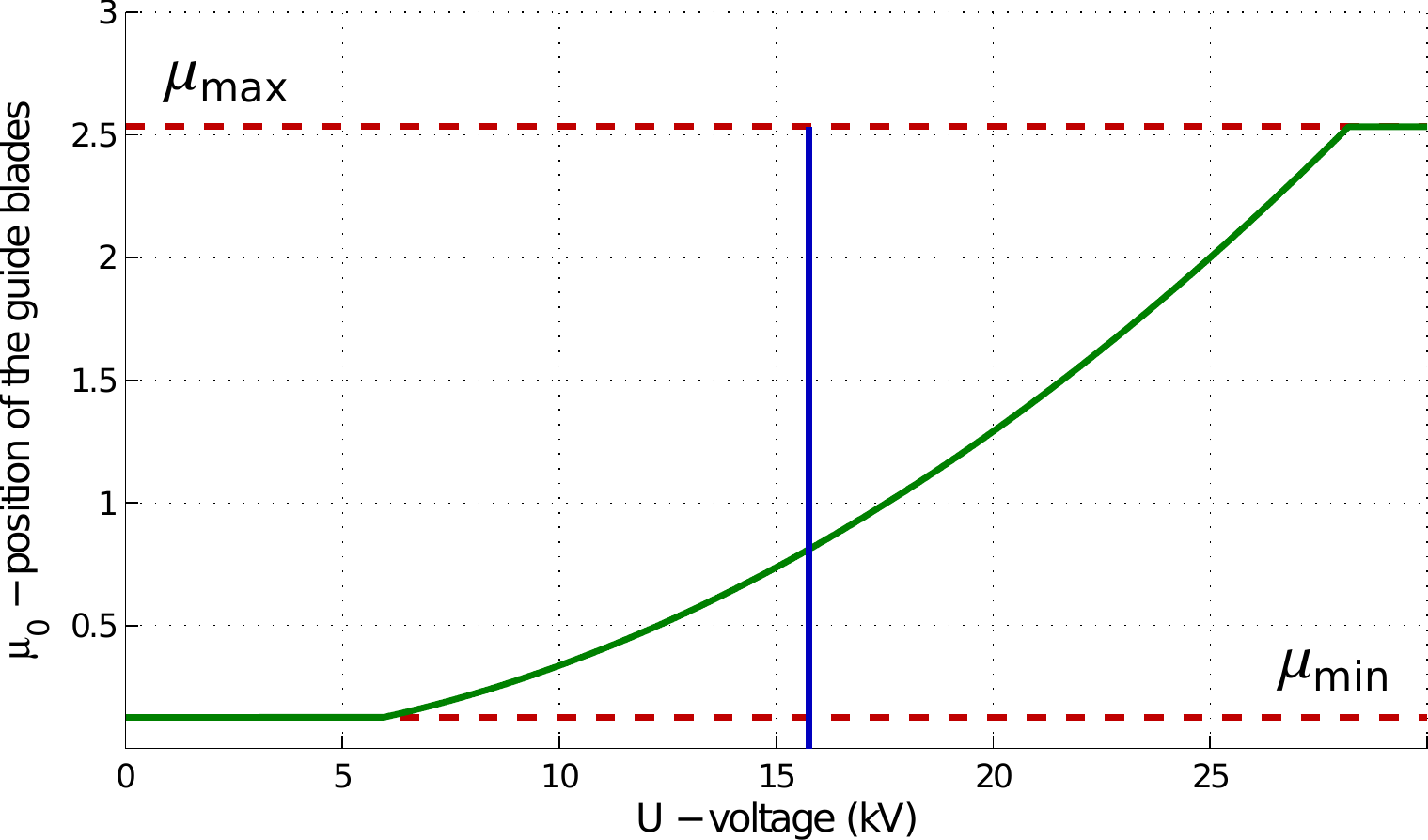}
	\caption{
		Plot of $\mu_0$ against $U$, $k=40$, $\theta_0 = \arccos(0.9)$
	}
	\label{mu_0U}
\end{figure}

Equation (\ref{balance}) contains the input parameter $U$ and the variable $\mu_0$. The equality of turbine and generator torques is
attained due to the variable  $\mu_0$. Recall that $\omega_0 = 14,954 \, \text{rad}/\text{s}$,  $\theta_0 = \arccos (0.9)$. Then $\mu_0$, depending on the voltage $U$, is found from the balance equation ~(\ref{balance}) of torques of turbine and generator. A plot of $\mu_0$ against $U$ is presented in Fig.~\ref{mu_0U}.

\begin{figure}[!ht]
	\begin{center}
		\includegraphics[width=1\linewidth]{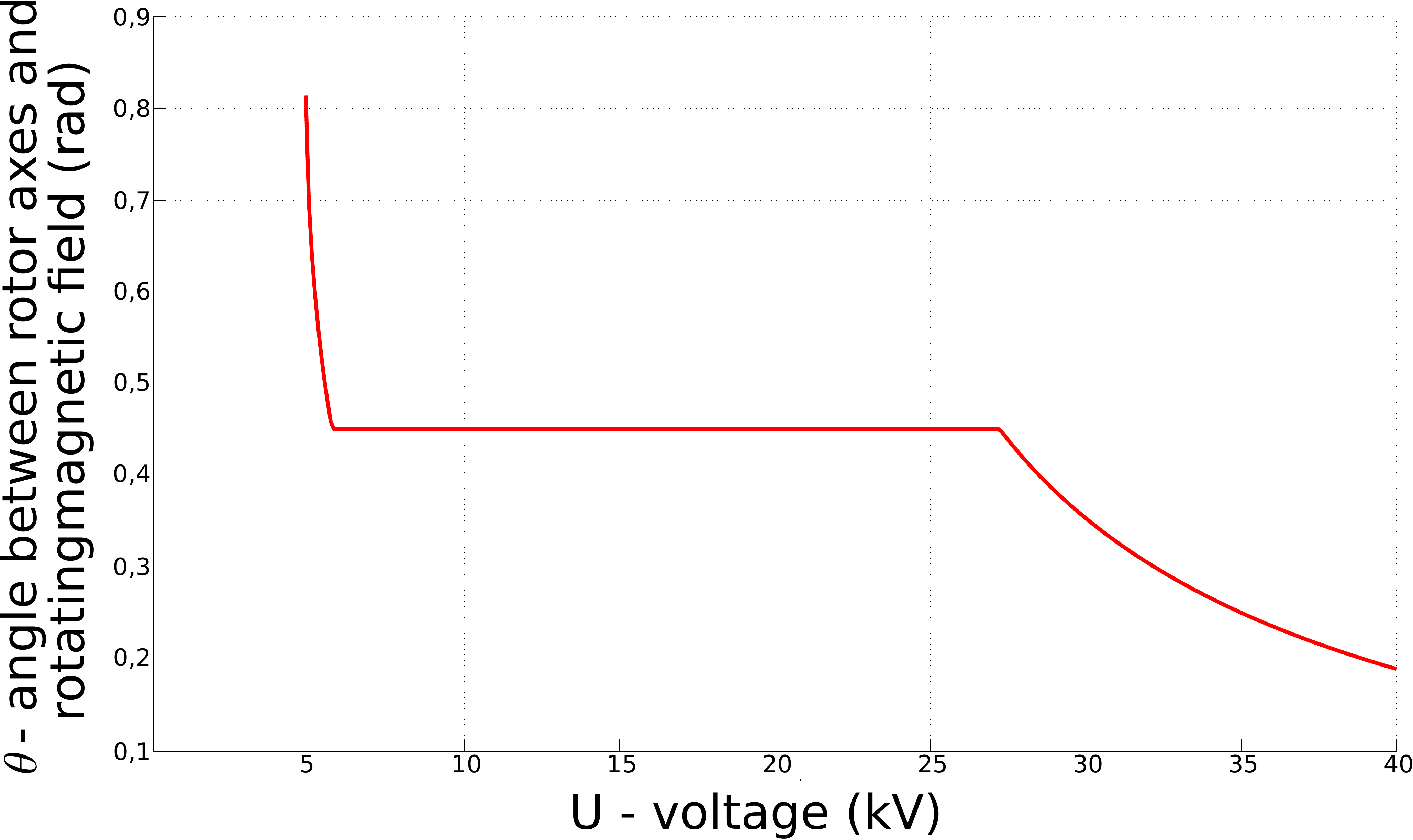}
	\caption{Plot of $\theta$ against $U$}
	\label{pic_static}
	\end{center}
\end{figure}

For sufficiently large and small values of $ U $, the parameter $ \mu_0 $ falls on the saturation that corresponds to the limit position of guide vanes. In particular this mode corresponds to the start of hydropower unit. 
Then it is necessary to consider a more complex balance equation
$$M_T(\mu_0 (U), \omega_0) = M_G(U,\theta_0 + \theta_\Delta),$$
and the balance will be achieved with help of $\theta_\Delta$, i.e. $\theta_\Delta^{st}$ may not be equal to zero. In Fig.~\ref{pic_static} a plot of $\theta=\theta_0+\theta_\Delta$ against voltage $U$ is shown.

If $\mu_0$ does not fall on the saturation, then $\theta_\Delta^{\text{st}} = 2\pi k$. Otherwise,
$\theta_\Delta^{\text{st}} = C + 2\pi k$, where  $C=\rm{const}$.

Since for $\mu_0 = \mu_{max}$ ($\mu_0 = \mu_{min}$) there are restrictions on $\theta$,
one can get again the saturation.
In this case the balance will be achieved with the help of $\omega_\Delta$.

However for the considered allowed voltage $U$ the value $\mu_0(U)$ does not fall on the saturation,
thus we do not consider the cases of achieving the balance with the help of  $\theta_\Delta$ and $\omega_\Delta$.

The instantaneous power is determined by the formula
$$
	P(U, \theta(t)) = - \frac{3}{2} \left( i_d(t) U \sin(\theta(t)) + i_q(t) U \cos(\theta(t)) \right).
$$
Represent the instantaneous power in the following form
$$
P(U, \theta(t)) = P_0  + P_\Delta(U, \theta(t)),
$$
where $P_0$ is the required (nominal) power, which is determined by the formula  $P_0=P(U, \theta_0)$, $P_\Delta$ is a deviation of the
required power.

\begin{figure}[H] \label{PU}
	\centering
	\includegraphics[width=0.9\linewidth]{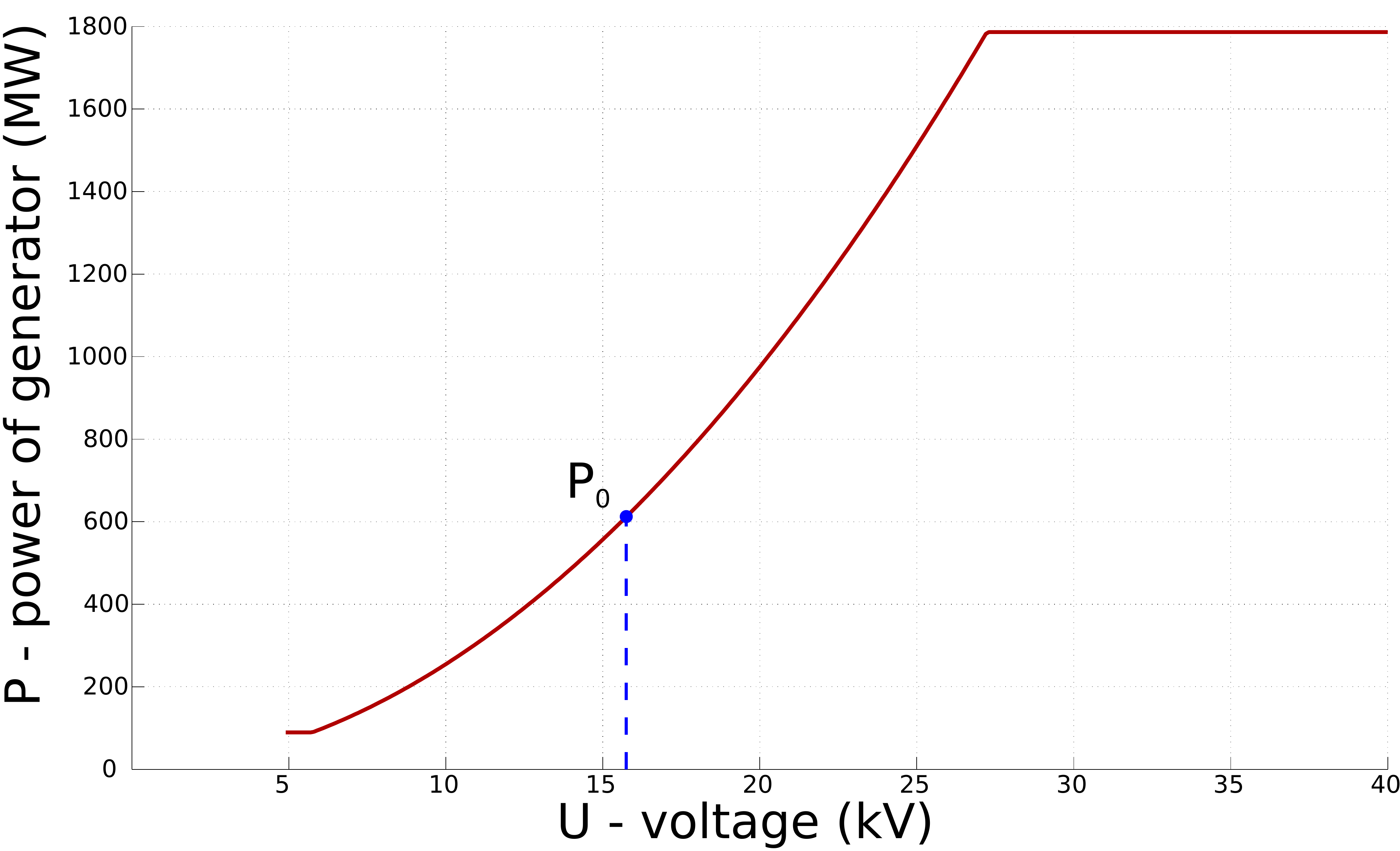}
	\caption{
		A plot of $P(U,\theta_0)$
	}
\end{figure}

The change of $U$ leads to the change of power. A plot of $P$ against $U$ is shown in
Fig.~\ref{PU}.

\section{Calculation of generator and turbine parameters for the Sayano-Shushenskaya hydropower plant}
\label{parameters}

The radial-axial vertical hydraulic turbines RO-230/833-B-677, connected
with synchronous generator on the umbrella type
SVF-1285/275-42 UHL4, are installed at the Sayano-Shushenskaya hydropower plant.

\begin{table}[h]
\caption{Parameters of synchronous generators (SVF-1285/275-42 UHL4) \cite{Karapetyan}}
	\begin{tabular}{|l|c|}
	\hline
		Parameter & Value
	\\ \hline
		Rated angular & $142.8$ [rad/s]
	\\  speed, $\omega_0$ &
	\\ \hline
		Stator resistance, $r$ & $0.0034$ [p.u.]
	\\ \hline
		Leakage inductive reactance & $0.184$ [[p.u.]
	\\  of stator winding, $x_s$ &
	\\ \hline
		Synchronous inductance & $1.58$ [p.u.]
	\\  along the axis $d$, $x_d$ &
	\\ \hline
		Synchronous inductance  & $0.97$ [p.u.]
	\\  along the axis $q$, $x_q$ &
	\\ \hline
		Transient resistance, $x_d'$ & $0.43$ [p.u.]
	\\ \hline
		Sub-transient reactance & $0.3$ [p.u.]
	\\   along the axis $d$, $x_d''$ &
	\\ \hline
		Sub-transient reactance & $0.31$ [p.u.]
	\\  along the axis $q$, $x_q''$ &
	\\ \hline
		Time constant & $8.21$ [s]
	\\  of field winding, $T_r$ &
	\\ \hline
		Sub-transient time & $0.143$ [s]
	\\  of field winding, $T_d''$ &
	\\ \hline
		Moment of inertia, $J$ & $25.5 \cdot 10^6$\,[$\text{kg}\cdot \text{m}^2$]
	\\ \hline
		Field voltage, $E_r$ & $530$ [V]
	\\ \hline
		Stator voltage, $U$ & $15.75$ [kV]
	\\ \hline
	\end{tabular}
\end{table}

The rest of the system parameters are determined as follows:

\begin{itemize}
	\item Impedances of stator winding along the axis $d$:
	$ x_{ad} = x_d - x_s = 1.396,
	$
	\item Impedances of stator winding along the axis $q$:
	$ x_{aq} = x_q - x_s = 0.786,
	$
	\item Impedances of field winding:
	$x_r = \frac{x_{ad}^2}{x_d - x_d'} = 1.6946,
	$
	\item Impedances of damper winding along the axis $d$:
	$
		x_{rd} = x_{ad} +
			\left(
				\frac{1}{x''_d - x_s} - \frac{1}{x_{ad}} - \frac{1}{x_{sr}}
			\right)^{-1} = 1.6155,
	$
	\item Impedances of damper winding along the axis $q$:
	$
		x_{rq} = x_{aq} +
			\left(
				\frac{1}{x''_q - x_s} - \frac{1}{x_{aq}}
			\right)^{-1} =  0.9361,
	$
	\item Resistance of field winding:
	$
		x_{sr} = x_r - x_{ad} = 0.2986,
	$
	\item Resistance of damper winding along the axis $d$:
	$
		r_{rd} = \frac{(x_{rd}x_d-x_{ad}^2)x_{rd}}{\omega_0 x_d x'_d T''_d} = 0.1246,
	$
	\item Resistance of damper winding along the axis $q$:
	$
		r_{rq} = \frac{x_{rq}x_q - x_{aq}^2}{\omega_0 x_q T''_q} = 0.0823,
	$
	\item Time constant of damper winding along the axis $d$:
	$
		T_{rd} = \frac{x_{rd}}{\omega_0 r_{rd}} = 0.8666,
	$
	\item Time constant of damper winding along the axis $q$:
	$
		T_{rq} = \frac{x_{rq}}{\omega_0 r_{rq}} = 0.7604.
	$
\end{itemize}

\begin{table}[h]
\caption{Parameters of turbine (RO-230/833-B-677)}
	\begin{tabular}{|l|c|}
	\hline
		Parameter & Value
	\\ \hline
		Length of penstock, $l$ & $212$ [m]
	\\ \hline
		Cross sectional diameter of penstock, $D$ &  $7.5$ [m]
	\\ \hline
		Available water flow through turbine, $Q_{nom}$  & $358$ [${\text{m}^3}/{\text{s}}$]
	\\ \hline
	\end{tabular}
\end{table}

Thus, for modeling the parameters of hydropower unit of the Sayano-Shushenskaya hydropower plant were used (\cite{Karapetyan}):
$\omega_0 = 2\pi 142.8/ 60$ [rad/s],
$r = 0.0034$ [p.u.],
$x_d = 1.58$ [p.u.],
$x_q = 0.97$ [p.u.],
$T_r = 8.21$ [s],		
$J = 25.5 \cdot 10^6$ [$\text{kg}\cdot \text{m}^2$],
$E_r = 530$ [p.u.],
$C = 0.27$ [$\text{m}^3\sqrt{\text{m}}/ \sqrt{\text{kg}}]$,
$x_{ad} = 1.396$ [p.u.],
$x_{aq} = 0.786$ [p.u.],
$x_r    = 1.6946$ [p.u.],
$x_{rd} = 1.6155$ [p.u.],
$x_{rq} = 0.9361$ [p.u.],
$T_{rd} = 0.8666$ [s],
$T_{rq} = 0.7604$ [s],
$S = \pi / 4 \cdot 7.5^2$ [$\text{m}^2$],
$l = 192$ [m],
$\rho = 0.98 \cdot 10^3 [\text{kg/m}^3]$,
$p_u = 2.7 \cdot 10^6$ [Pa],
$p_l = 0.35 \cdot 10^6$ [Pa],
$k = 40$ [$\text{kg}\cdot \text{m}^2/ \text{s}^2$],
$Q_\text{max} = 358 \, [\text{m}^3/\text{s}]$.

\section{Local analysis}

Let us study the local stability of equilibrium points of system (\ref{model}). 
It is enough to carry out the analysis of stability on the interval $[0, 2\pi)$
since the solutions of the system are $ 2\pi $-periodic. The equilibrium points with respect to $\theta$ are defined from balance equation (\ref{balance}). 
System (\ref{model}) may have 0, 1, 2, 3 or 4 equilibrium points on the interval $ [0, 2\pi)$.

Let us find the Jacobian matrix of the right-hand side of system (\ref{model}).
For this purpose algebraic system of equations (\ref{equation_flux}) is solved for  $i_d$, $i_q$, $E_q$, $E_{rd}$, $E_{rq}$:
$$i_d    = X_d \Psi_d - X_r \Psi_r + X_{rd} \Psi_{rd}, \quad i_q    = Y_q \Psi_q - Y_{rq} \Psi_{rq},$$
$$E_q    = Z_d \Psi_d + Z_r \Psi_r - Z_{rd}, \Psi_{rd}, \quad E_{rd} = P_q \Psi_q - P_{rq} \Psi_{rq}, $$
$$E_{rq} = - Q_d \Psi_d + Q_r \Psi_r + Q_{rd} \Psi_{rd}, $$
where
$$X_{d}  = \frac{a_{4}\,a_{6} - 1}{b_{1}}, \quad
X_{r}  = \frac{a_{6} - 1}{b_{1}}, \quad
X_{rd} = \frac{1 - a_{4}}{b_{1}}, $$
$$Y_{q}  = \frac{1}{a_{2} - a_{7}}, \quad
Y_{rq} = \frac{1}{a_{2} - a_{7}},$$
$$Z_{d}  = \frac{(b_{1} - (a_{1}-a_{5})\,(a_{4}\,a_{6} - 1) )}{b_{1}\,(1 - a_{6})},\quad
Z_{r}  = \frac{a_{5} - a_{1}}{b_{1}},$$
$$Z_{rd} = \frac{b_{1} + (a_{1}-a_{5})\,(1 - a_{4})}{b_{1}\,(1 - a_{6})},$$
$$P_{q}  = \frac{a_{7}}{a_{2} - a_{7}}, \quad
P_{rq} = \frac{a_{2}}{a_{2} - a_{7}},$$
$$Q_{d}  = \frac{(a_{1} - a_{3})\,(a_{4}\,a_{6} - 1) - b_{1}}{(b_{1}\,(1 - a_{4})},$$
$$Q_{r}  = \frac{(a_{1} - a_{3})\,(a_{6} - 1) - b_{1}}{b_{1}\,(1 - a_{4})}, \quad
Q_{rd} = \frac{a_{3} - a_{1}}{b_{1}},$$
$$a_{1} = x_{d}, \quad
a_{2} = x_{q}, \quad
a_{3} = \frac{x_{ad}^2}{x_{r}}, \quad
a_{4} = \frac{x_{ad}}{x_{r}},$$
$$a_{5} = \frac{x_{ad}^2}{x_{rd}}, \quad
a_{6} = \frac{x_{ad}}{x_{rd}}, \quad
a_{7} = \frac{x_{aq}^2}{x_{rq}},$$
$$b_{1} = (a_{1} - a_{5})\,(a_{4} - 1) - (a_{1}\,a_{4} - a_{3})\,(1 - a_{6}).$$

Then nonzero elements of the Jacobi matrix
$$J = \{j_{i,k}\}^{i=1...9}_{k=1...9}$$
of the right-hand side of system (\ref{model}) in stationary point are defined by the formulas:
$$j_{1,2} = \omega_0, \quad
j_{2,2} = \frac{-k\,(Q^{st})^3}{T_J\,C^2\,\omega_0^2\,\mu_0 ^2}, \quad
j_{2,3} = \frac{3\,k\,(Q^{st})^2}{T_J\,C^2\,\omega_0^2 \,\mu_0^2},$$
$$j_{2,4} = \frac{(X_d - Y_q)\,\Psi_q + Y_{rq}\,\Psi_{rq}}{T_J},$$
$$j_{2,5} = \frac{-Y_q\,\Psi_d + X_d\,\Psi_d - X_r\,\Psi_r + X_{rd}\,\Psi_{rd}}{T_J},$$
$$j_{2,6} = \frac{-X_r\,\Psi_q}{T_J}, \quad
j_{2,7} = \frac{X_{rd}\,\Psi_q}{T_J}, \quad
j_{2,8} = \frac{Y_{rq}\,\Psi_d}{T_J},$$
$$j_{2,9} = \frac{2\,k\,(Q^{st})^3}{T_J\,C^2\,\mu_0^3\,\omega_0^2},\quad
j_{3,3} = -\frac{2\,S\,Q^{st}}{l\,\rho\,C^2\,\mu_0^2},$$
$$j_{3,9} = -\frac{2\,S\,(Q^{st})^2}{l\,\rho\,C^2\,\mu_0^3}, \quad
j_{4,1} =  \omega_0\,U\,\cos(\theta_0+\theta_\Delta^{st}),$$
$$j_{4,2} = -\omega_0\,\Psi_q, \quad
j_{4,4} = -\omega_0\,r\,X_d, \quad
j_{4,5} = -\omega_0$$
$$j_{4,6} =  \omega_0\,r\,X_r, \quad
j_{4,7} = -\omega_0\,r\,X_{rd},$$
$$j_{5,1} =  \omega_0\,U\,\sin(\theta_0+\theta_\Delta^{st}), \quad
j_{5,2} =  \omega_0\,\Psi_d,$$
$$j_{5,4} =  \omega_0, \quad
j_{5,5} = -\omega_0\,r\,Y_q, \quad
j_{5,8} =  \omega_0\,r\,Y_{rq},$$
$$j_{6,4} = -\frac{Z_d}{T_r},\quad
j_{6,6} = -\frac{Z_r}{T_r},\quad
j_{6,7} =  \frac{Z_{rd}}{T_r},$$
$$j_{7,4} =  \frac{Q_d}{T_{rd}},\quad
j_{7,6} = -\frac{Q_r}{T_{rd}}, \quad
j_{7,7} = -\frac{Q_{rd}}{T_{rd}},$$
$$j_{8,5} =  \frac{P_q}{T_{rq}},\quad
j_{8,8} = -\frac{P_{rq}}{T_{rq}}, \quad
j_{9,2} =  \frac{1}{T_c}, \quad
j_{9,9} = -\frac{1}{T_c}.$$

\begin{figure}[!ht]
	\begin{center}
		\includegraphics[width=1\linewidth]{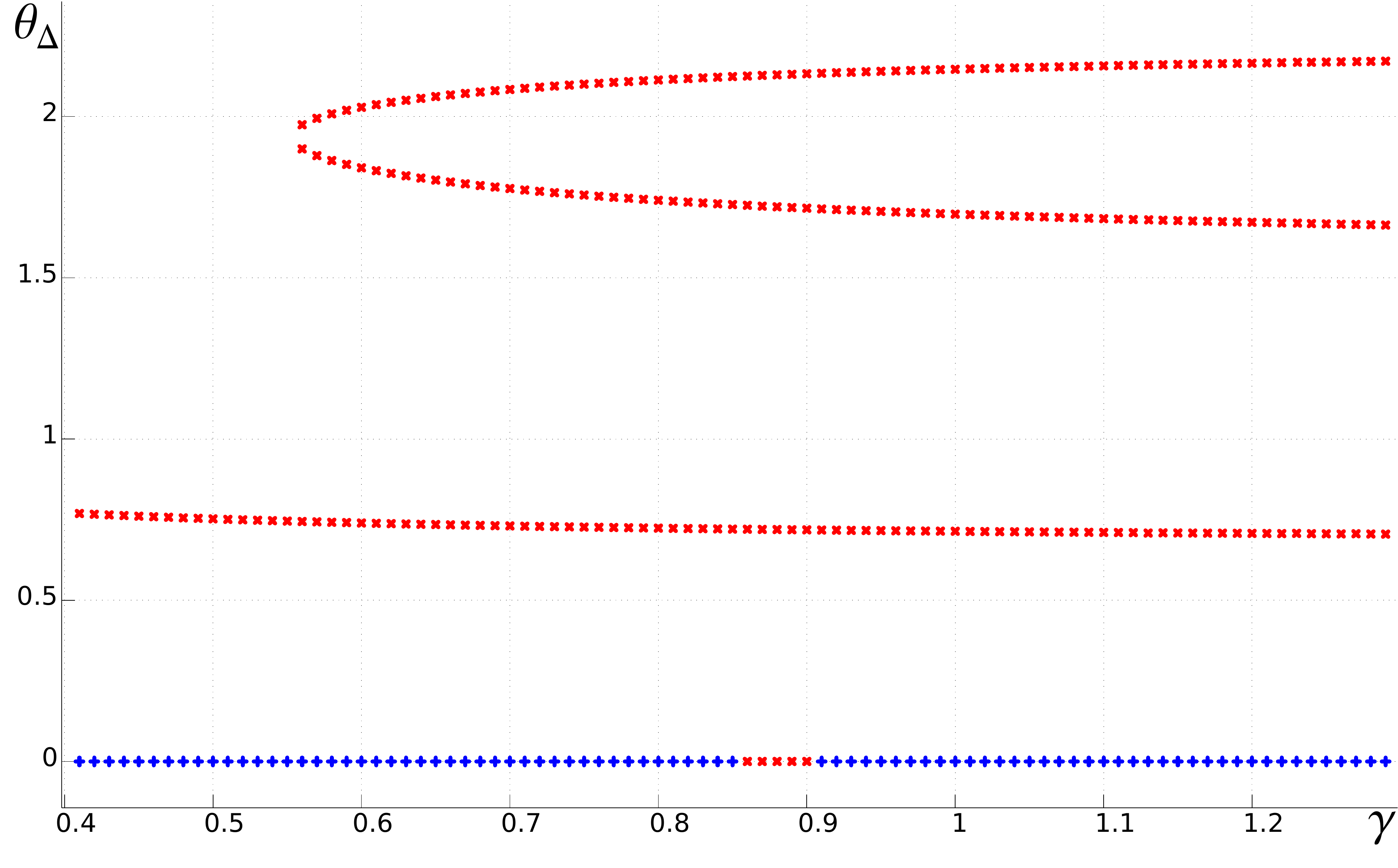}
	\caption{Stability of equilibrium states of system (\ref{model}), defined by the least-squares method: blue pluses are stable equilibrium states, red crosses are unstable equilibrium states}
	\label{eq_lsqnonlin}
	\end{center}
\end{figure}

In order to study the local stability of the equilibrium states of system (\ref{model}) relative to the voltage, the standard function 
\emph{lsqnonlin} of the application package {\sc MatLab} was used\footnote{
First, the study of local stability of equilibria was carried out 
via eigenvalues of the Jacobian matrix and the standard function \emph{eig} from {\sc MatLab}. 
It turned out that all equilibrium states are unstable for all $\gamma$. 
Further, the stability criterion of  Routh--Hurwitz for the Jacobian matrix was applied. It showed that the equilibrium state $\theta = \theta_0$ is unstable for $\gamma \in [\gamma_1 \approx 0.84; \gamma_2 \approx 1.17]$ and stable for other $\gamma$.
The obtained results of these methods can be explained by the fact that all calculations are made with some numerical error 
and because of the substantial difference in magnitude.} 
This function is based on the least-squares method, i.e. on  an iterative approximation to equilibrium state, that allows one to reduce the computational error and to define more exactly the interval of instability of the system. The initial data for this method is putative equilibrium state.
Using this method the Jacobi matrix calculated at an equilibrium is found,
and then the Routh--Hurwitz stability criterion is applied.
Results are presented in Fig.~\ref{eq_lsqnonlin}: the equilibrium state $ \theta = \theta_0 $ is unstable for $ \gamma \in [\gamma_1 \approx 0.86; \gamma_2 \approx 0.9] $ and stable for other $ \gamma $, the rest three equilibria states are always unstable.

\section{Analysis of transient processes}

During the operation of hydropower unit the transient processes related to sudden changes of the work parameters of the hydraulic unit often occur. As a result, the following problem arises: to find parameters, under which the hydropower unit pulls in the new operating mode after transient processes. This problem is closely related to the limit (ultimate) load problem, which arises in practice of operation of electrical motors after sudden change of load torque on the shaft \cite{1950-Annett,1958-Yanko-Trinitskii,1995-Haque,2002-Das,2005-Bianchi,LeonovK-2009}. For its solution the equal-area method is widely used in engineering practice. 
This method was used for some models in the works of A.A. Yanko-Trinitskii \cite{1958-Yanko-Trinitskii}. 
In our work modern methods of numerical integration of the system (Runge--Kutta method) are combined
with the analysis in the spirit of the classical ideas of Yanko--Trinitskii.

First the numerical analysis of transient processes was carried out with the initial data taken from a small neighborhood of the equilibrium state. It is verified whether the trajectory goes out from this neighborhood after a long integration time (1000 s) or not. As a result of the study it was obtained that the equilibrium state $ \theta = \theta_0 $, corresponding to the operating mode, is unstable for $ \gamma \in [\gamma_1 \approx 0.85; \gamma_2 \approx 0.91]$. The received interval of instability is consistent with the local analysis and corresponds to the interval [$ S_1 $, $ S_2 $] in Fig. ~\ref{zones_1}.
It can be defined more exactly due to coefficient $ k $, corresponding the turbine used at the Sayano-Shushenskaya hydropower plant.

Further, the numerical analysis of transient processes was carried out with various initial data. 
For values $\gamma$, corresponding the local stability, hidden oscillations\footnote{
 An oscillation can generally be easily numerically localized if the initial data from its open
 neighborhood in the phase space lead to
 a long-term behavior that approaches the oscillation.
 Therefore, from a computational perspective, it is natural to suggest the following classification
 of attractors \cite{KuznetsovLV-2010-IFAC,LeonovKV-2011-PLA,LeonovKV-2012-PhysD,LeonovK-2013-IJBC,KuznetsovL-2014-IFACWC},
 which is based on the simplicity of finding their basins of attraction in the phase space:
 {\it An attractor is called a \emph{self-excited attractor}
 if its basin of attraction
 intersects with any open neighborhood of an equilibrium,
 otherwise it is called a \emph{hidden attractor}.
 }
 For a \emph{self-excited attractor} its basin of attraction
 is connected with an unstable equilibrium
 and, therefore,
  (\emph{standard computational procedure})
 self-excited attractors
 can be localized numerically by the
\emph{standard computational procedure}:
 by constructing a solution using initial data from
 an unstable manifold in a neighborhood of an unstable equilibrium,
 and observing how it is attracted, and visualizing the oscillation.
 In contrast,  the basin of attraction for a hidden attractor is not connected with any equilibrium.
 For example, hidden attractors are attractors in systems
 with no equilibria or with only one stable equilibrium
 (a special case of the multistability: coexistence of attractors in multistable systems).
 Well known examples of the hidden oscillations are nested limit cycles in 16th Hilbert problem
(see, e.g., \cite{KuznetsovKL-2013-DEDS,LeonovK-2013-IJBC})
and counterexamples to the Aizerman and Kalman conjectures on the absolute stability of nonlinear control systems \cite{LeonovBK-2010-DAN,LeonovK-2011-DAN,BraginVKL-2011,LeonovK-2013-IJBC}. 
Hidden oscillations in the models of electrical machines are discussed, e.g., in \cite{KiselevaKLN-2012-IEEE,KiselevaKLN-2013,LeonovK-2013-IJBC,LeonovKKSZ-2014,KiselevaKKLS-2014}}
are not found numerically and all simulated trajectories attract to the equilibrium states. 
For values $\gamma$, corresponding to the local instability, simulated trajectories attract to self-excited periodic solutions. 
The local bifurcation, in which an equilibrium loses stability and a small stable limit cycle is born, occurs in considered multidimensional system (this bifurcation is an analog of Andronov-Hopf bifurcation \cite{andronov1971theory,guckenheimer1983nonlinear,kuznetsov2013elements}).

According to \cite{Act}, immediately before the accident the power of the second hydropower unit was 475~MW at a head of 212 meters 
(i.e., it worked in the not recommended zone II (Fig.~\ref{zones_1}).
On the day of the accident the power of the second hydropower unit was reduced in accordance with the commands of  
the group controller of active and reactive power.


Below three cases are modeled:
\begin{enumerate}
\item operation of hydropower unit at the rated voltage (Fig.~\ref{zones_1}, point~A) with initial data \\
$(\theta_\Delta, s, Q, \Psi_d, \Psi_q, \Psi_r, \Psi_{rd}, \Psi_{rq}, \mu_\Delta) = \\ $ $ \text{~~~~~~~~~~~~~~~~~~~~~~~~~~~~~~~~~~~~~~~~} (0, 1, 0, 0, 0, 0, 0, 0, 0)$,
\item reducing the power of hydropower unit, that corresponds to reducing voltage to 0.89 of the rated voltage (Fig.~\ref{zones_1}, point~B),
\item reducing power of hydropower unit, that corresponds to reducing voltage to 0.7 of the rated voltage (Fig.~\ref{zones_1}, point~C).
\end{enumerate}

\begin{figure}[H]
	\begin{center}
	\includegraphics[width=1\linewidth]{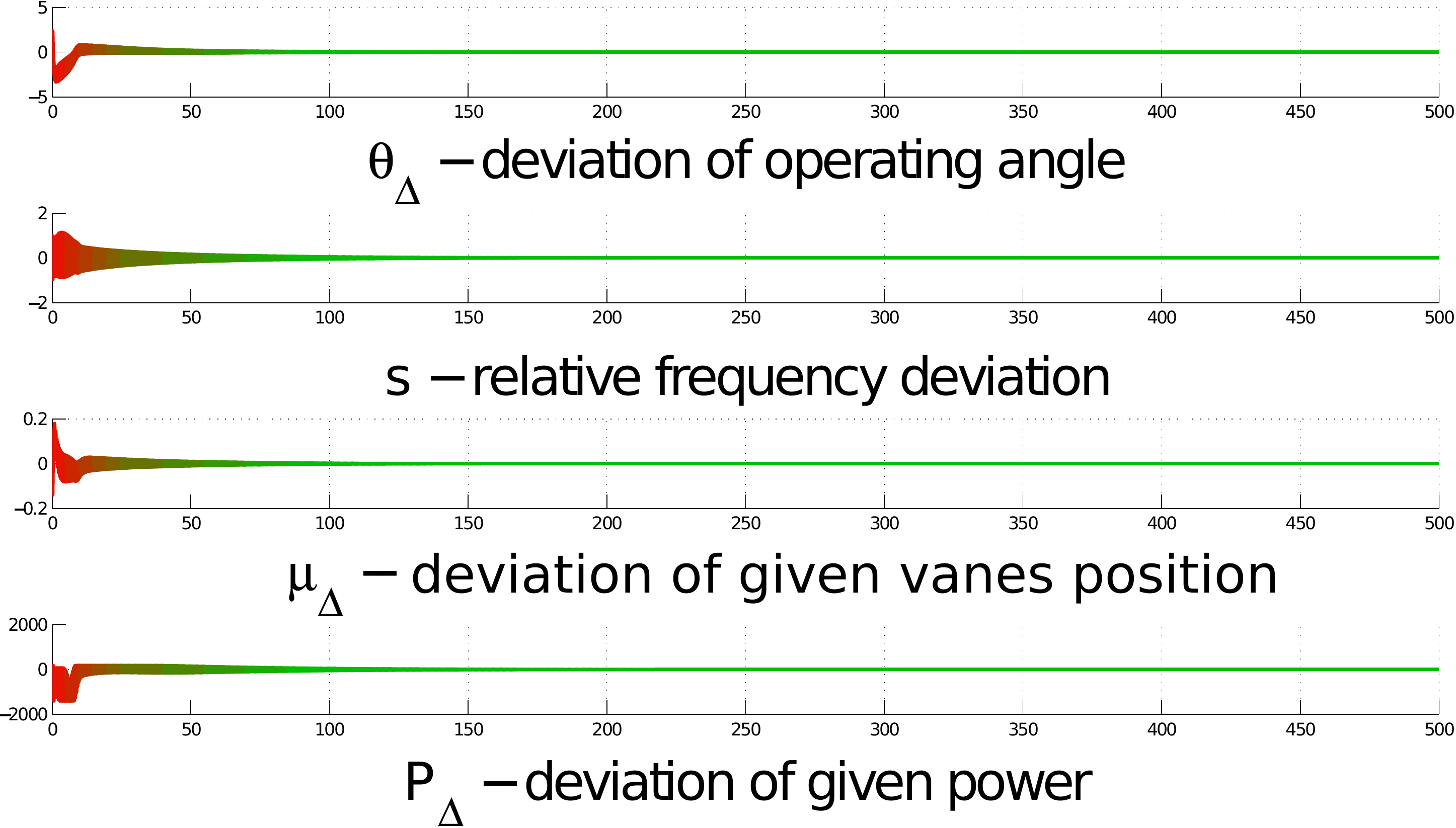}
\caption{Stable equilibrium in the mathematical model of hydropower unit, $U = 15.75 \cdot 10^3$ [V]} %
	\label{U2}	
	\end{center}
\end{figure}

\begin{figure}[H]
	\begin{center}
	\includegraphics[width=1\linewidth]{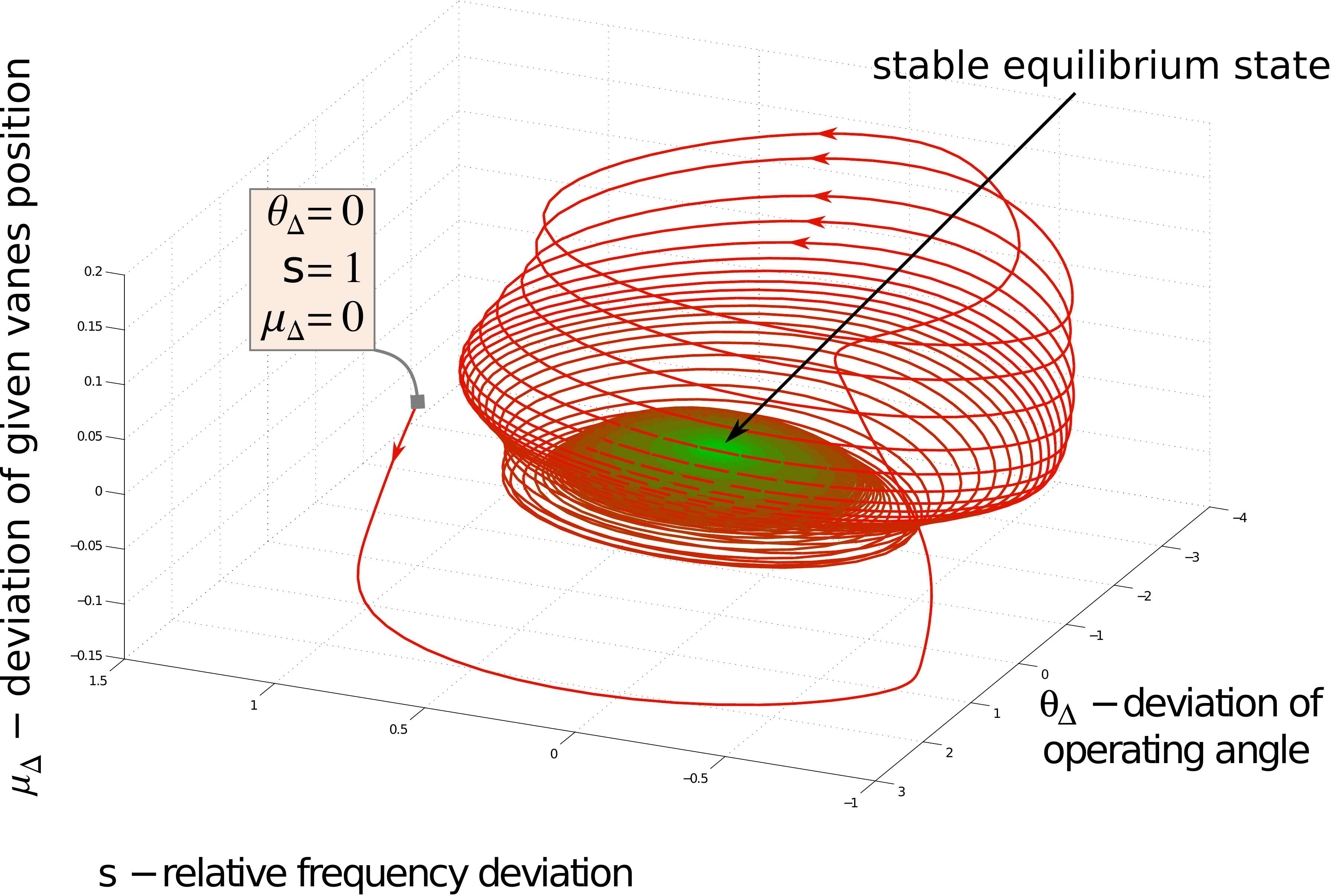}
\caption{Stable equilibrium in the mathematical model of hydropower unit -- projection onto ($\theta_\Delta, \mu_\Delta, s$), $U = 15.75 \cdot 10^3$ [V]}
	\label{U3}	
	\end{center}
\end{figure}

\begin{figure}[H]
	\begin{center}
	\includegraphics[width=1\linewidth]{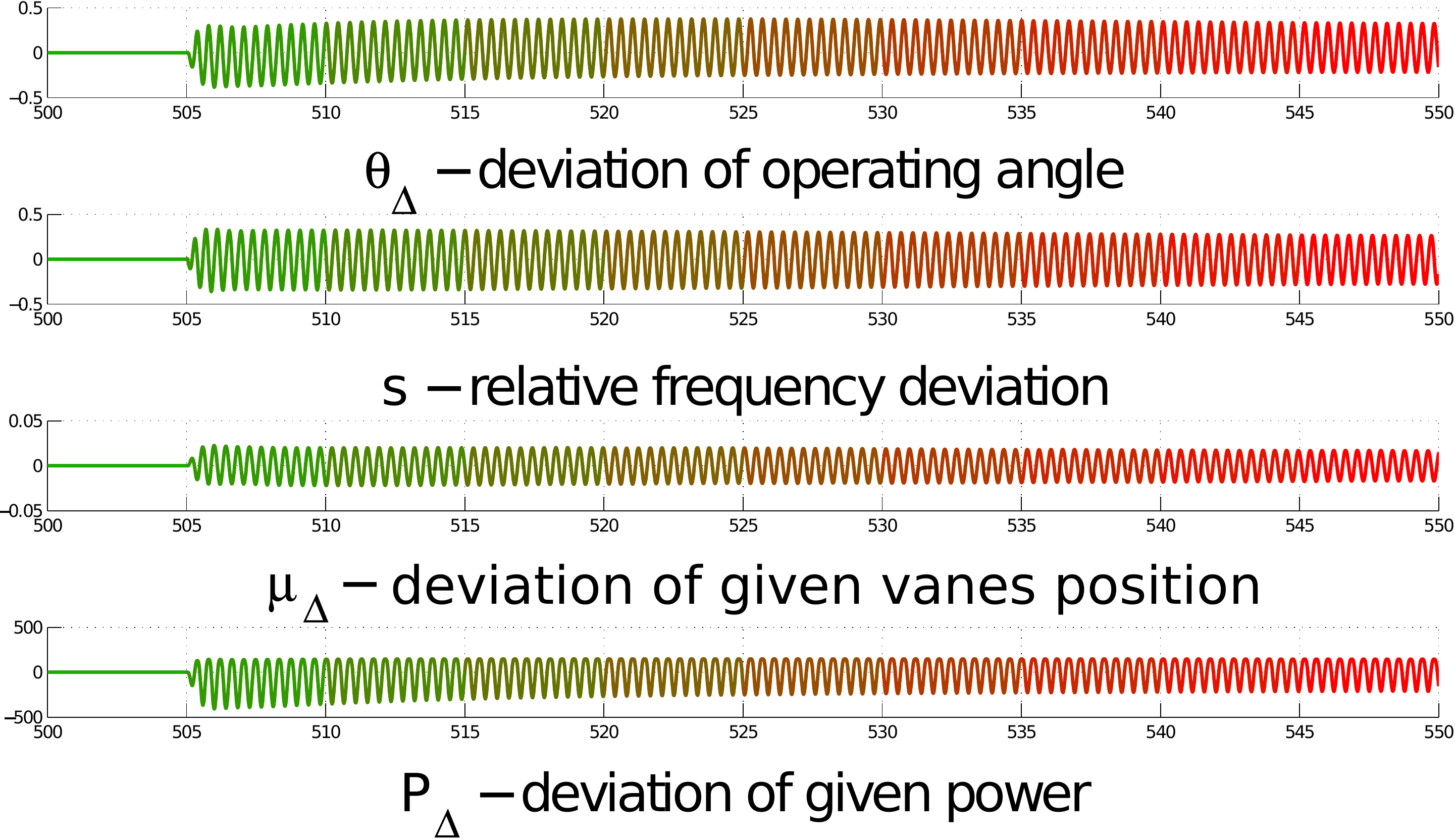}
\caption{Limit cycle in the mathematical model of hydropower unit, $U = 0.89 \cdot 15.75 \cdot 10^3$ [V]} %
	\label{0.9U2}	
	\end{center}
\end{figure}

\begin{figure}[H]
	\begin{center}
	\includegraphics[width=1\linewidth]{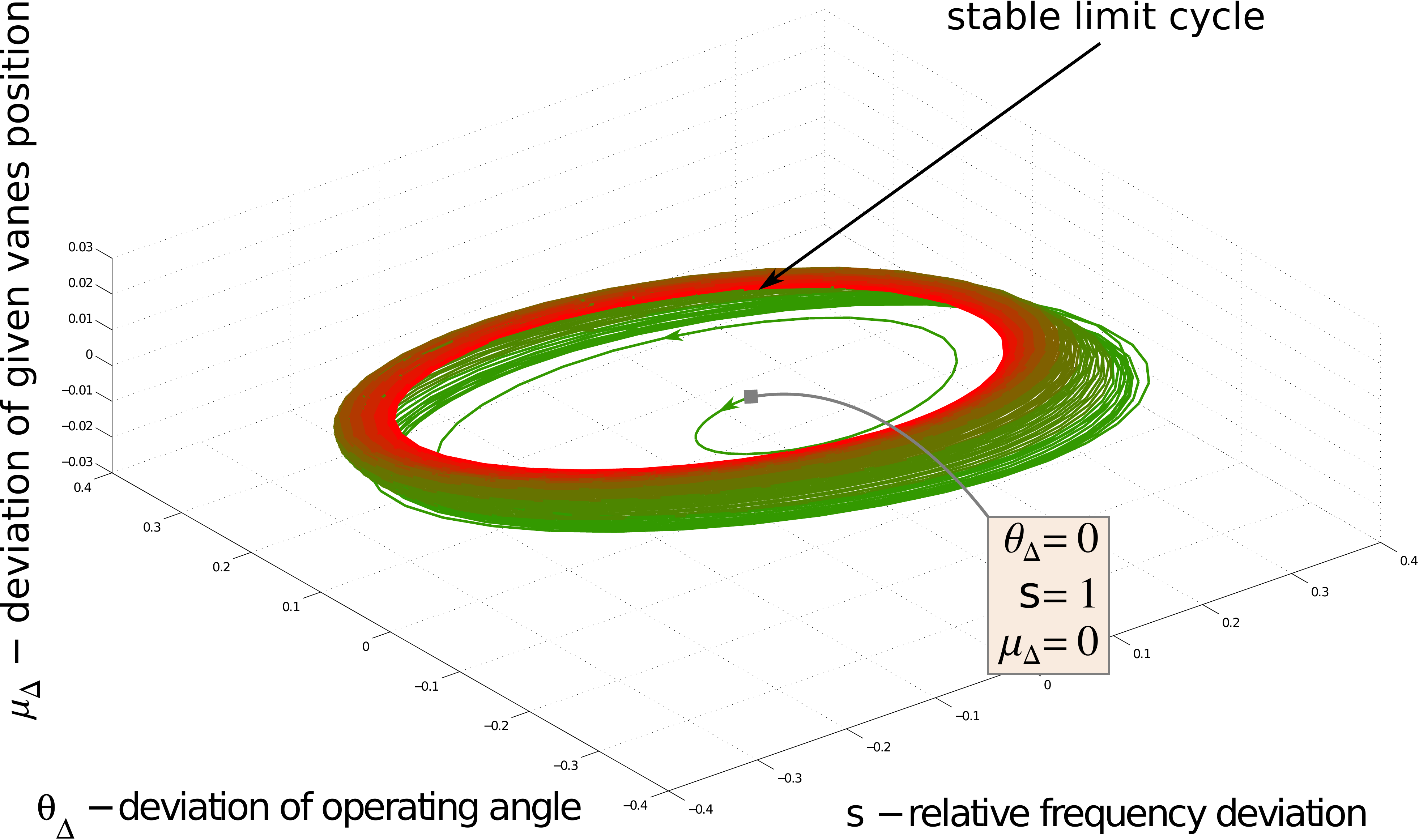}
\caption{Limit cycle in the mathematical model of hydropower unit -- projection onto ($\theta_\Delta, \mu_\Delta, s$), $U = 0.89 \cdot 15.75 \cdot 10^3$ [V]}
	\label{0.9U3}	
	\end{center}
\end{figure}

\begin{figure}[H]
	\begin{center}
	\includegraphics[width=1\linewidth]{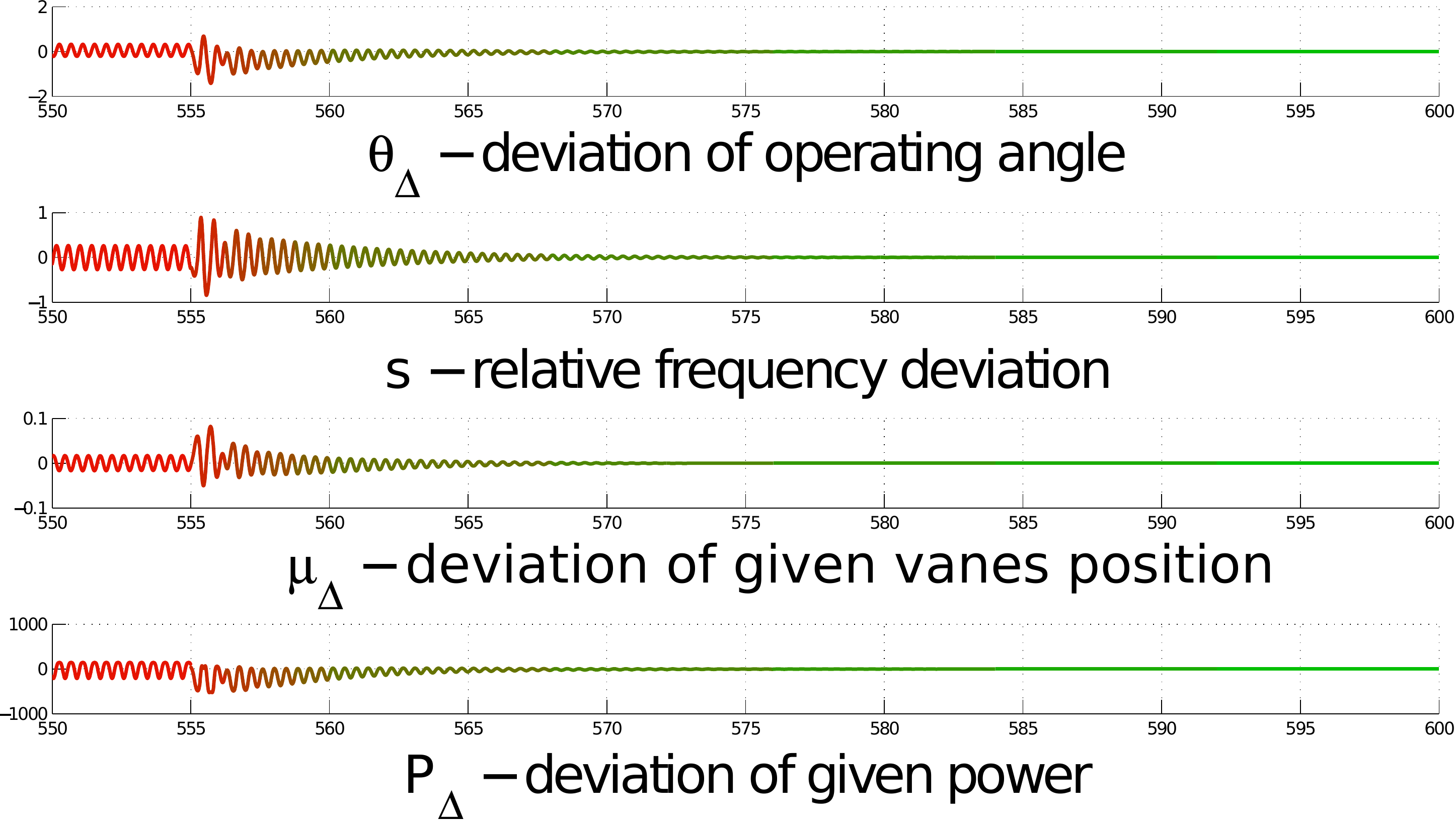}
\caption{Stable equilibrium in the mathematical model of hydropower unit, $U = 0.7 \cdot 15.75 \cdot 10^3$ [V]} %
	\label{0.8U2}	
	\end{center}
\end{figure}

\begin{figure}[H]
	\begin{center}
	\includegraphics[width=1\linewidth]{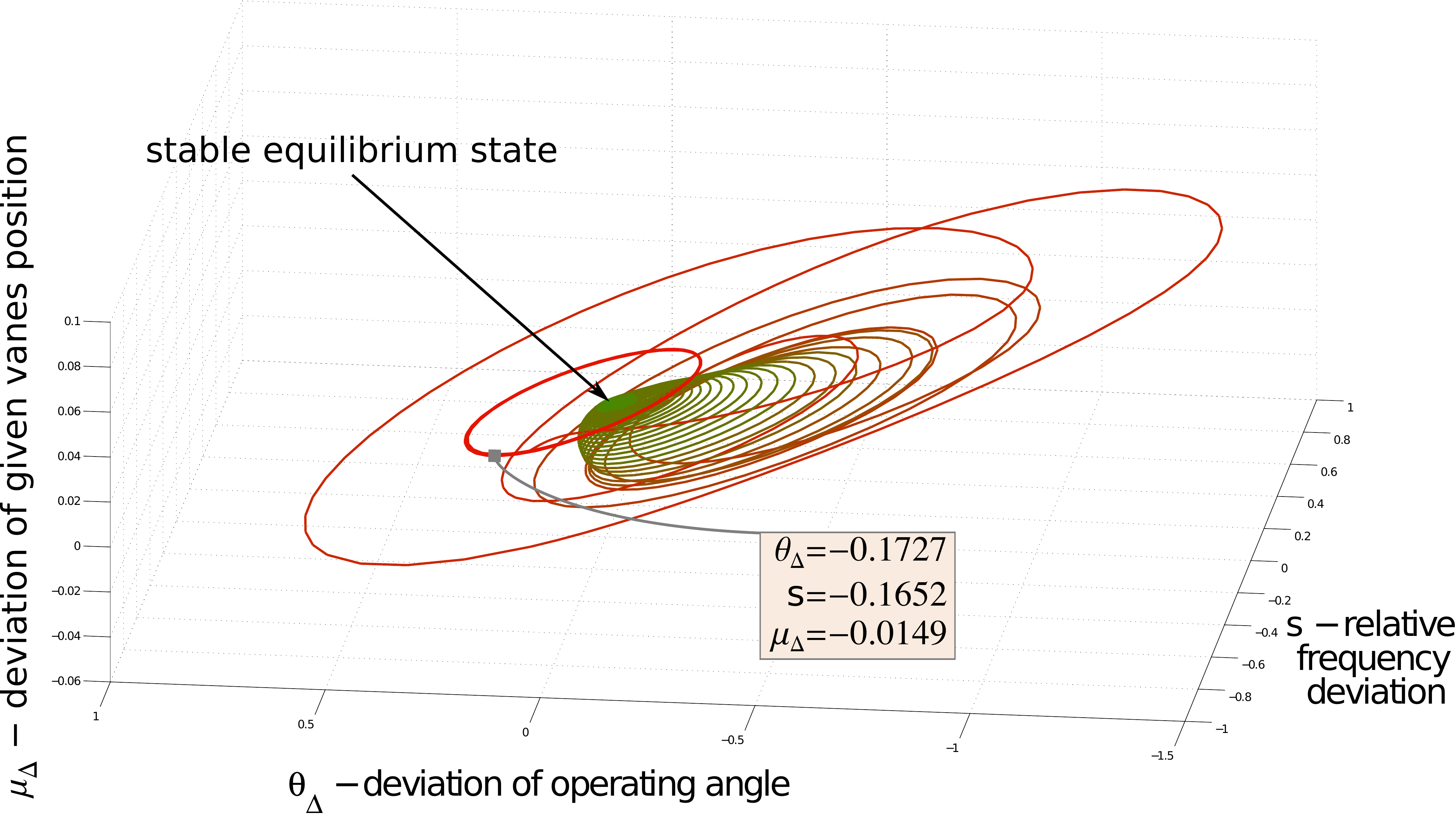}
\caption{Stable equilibrium in the mathematical model of hydropower unit -- projection onto ($\theta_\Delta, \mu_\Delta, s$), $U = 0.7 \cdot 15.75 \cdot 10^3$ [V]}
	\label{0.8U3}	
	\end{center}
\end{figure}

The results of modeling have shown that at the rated voltage the trajectory of the system after transient processes
is attracted to the equilibrium state, which corresponds to operating mode of hydropower unit (Figs. \ref{U2}, \ref{U3}).
Further at some instant the voltage is reduced to $0.89$ of the rated voltage.
In this case the trajectory of the system after transient processes is attracted to the stable limit cycle (Figs. \ref{0.9U2}, \ref{0.9U3}), i.e., vibrations are arisen in the hydropower units. Then the voltage is reduced to $0.7$ of the rated voltage. The trajectory of the system after transient processes is attracted to equilibrium state (Figs. \ref{0.8U2}, \ref{0.8U3}).

\begin{figure}[H] \label{amplituda}
	\centering
		\includegraphics[width=1\linewidth]{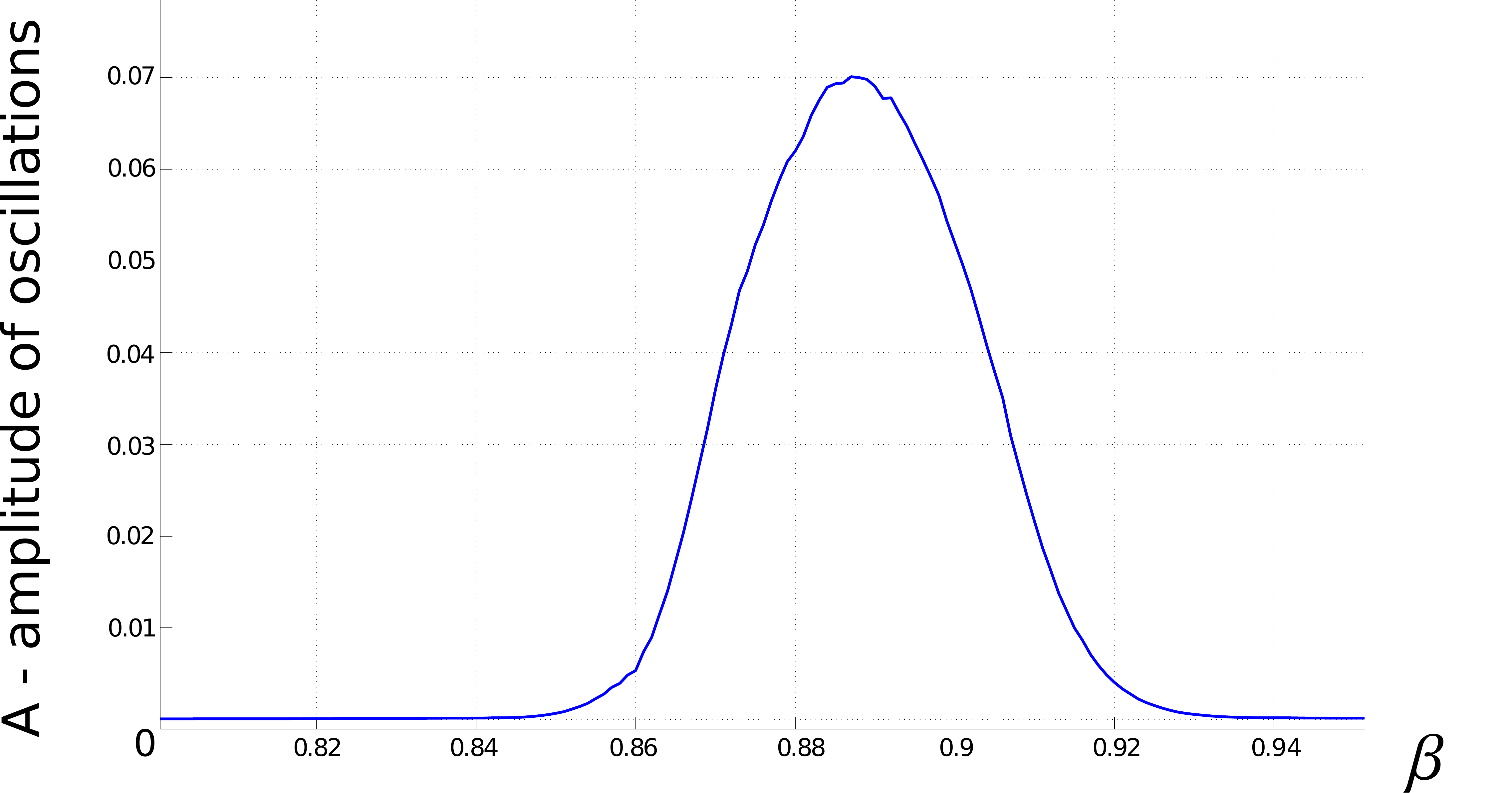}
		\caption{
		Amplitude of oscillations for different voltages $U = \beta U_ {nom}$ ($ \beta $ - percent of capacity)
		}
\end{figure}

From Fig.~\ref{amplituda}, it is clear that the maximum value of the oscillation amplitude is reached at $U = 0.89U_{nom}$.

Thus, the results of modeling are sufficiently consistent with the full-scale tests carried out for hydropower units of the Sayano-Shushenskaya hydropower plant.

\begin{acknowledgements}
The work was supported Russian Scientific Foundation
(14-21-00041) and Saint-Petersburg State University.
\end{acknowledgements}

\bibliographystyle{spmpsci}      

\begin{thebibliography}{10}
\providecommand{\url}[1]{{#1}}
\providecommand{\urlprefix}{URL }
\expandafter\ifx\csname urlstyle\endcsname\relax
  \providecommand{\doi}[1]{DOI~\discretionary{}{}{}#1}\else
  \providecommand{\doi}{DOI~\discretionary{}{}{}\begingroup
  \urlstyle{rm}\Url}\fi

\bibitem{Adkins}
Adkins, B.: The General Theory of Electrical Machines.
\newblock Chapman \& Hall, London (1962)

\bibitem{andronov1971theory}
Andronov, A.A.: Theory of bifurcations of dynamic systems on a plane, vol. 554.
\newblock Israel Program for Scientific Translations;[available from the US
  Dept. of Commerce, National Technical Information Service, Springfield, Va.]
  (1971)

\bibitem{1950-Annett}
Annett, F.A.: Electrical machinery: a practical study course on installation,
  operation and maintenance.
\newblock McGraw-Hill (1950)

\bibitem{2005-Bianchi}
Bianchi, N.: Electrical machine analysis using finite elements.
\newblock CRC Press (2005)

\bibitem{boldea}
Boldea, I.: Synchronous generators.
\newblock CRC Press (2006)

\bibitem{bondareva2004}
Bondareva, N., Kolotovkin, D., Cherkaoui, R., Germond, A., Grobovoy, A.,
  Stubbe, M.: Comparison of the results of full-scale experiment and long term
  dynamics simulation in the siberian interconnected power system.
\newblock In: Proc. IREP Symposium, Bulk Power System Dynamics and Control-VI,
  pp. 492--498 (2004)

\bibitem{BraginVKL-2011}
Bragin, V.O., Vagaitsev, V.I., Kuznetsov, N.V., Leonov, G.A.: Algorithms for
  finding hidden oscillations in nonlinear systems. {T}he {A}izerman and
  {K}alman conjectures and {C}hua's circuits.
\newblock Journal of Computer and Systems Sciences International
  \textbf{50}(4), 511--543 (2011).
\newblock \doi{10.1134/S106423071104006X}

\bibitem{2002-Das}
Das, J.C.: Power system analysis: short-circuit load flow and harmonics.
\newblock CRC Press (2002)

\bibitem{GaleaniOnoriTeel_SCL08}
Galeani, S., Onori, S., Teel, A.R., Zaccarian, L.: A magnitude and rate
  saturation model and its use in the solution of a static anti-windup problem.
\newblock Systems \& Control letters \textbf{57}(1), 1--9 (2008)

\bibitem{Gluhih}
Gluhih, V.: Full-scale tests of turbines of the {S}ayano-{S}hushenskaya hydro
  power plant with standard runners.
\newblock Tech. Rep. No. 1008 (1988)

\bibitem{guckenheimer1983nonlinear}
Guckenheimer, J., Holmes, P.: Nonlinear oscillations, dynamical systems, and
  bifurcations of vector fields, vol.~42.
\newblock Springer Science \& Business Media (1983)

\bibitem{1995-Haque}
Haque, M.H.: Further developments of the equal-area criterion for multimachine
  power systems \textbf{33(3)}, 175--183 (1995)

\bibitem{Karapetyan}
Karapetyan, I., Faibisovich, D., Sharipo, I.: Handbook for designing electrical
  networks.
\newblock Litres (2013)

\bibitem{KiselevaKKLS-2014}
Kiseleva, M., Kondratyeva, N., Kuznetsov, N., Leonov, G., Solovyeva, E.: Hidden
  periodic oscillations in drilling system driven by induction motor.
\newblock IFAC Proceedings Volumes (IFAC-PapersOnline) \textbf{19}, 5872--5877
  (2014)

\bibitem{KiselevaKLN-2013}
Kiseleva, M., Kuznetsov, N., Leonov, G., Neittaanmaki, P.: Hidden oscillations
  in drilling system actuated by induction motor.
\newblock IFAC Proceedings Volumes (IFAC-PapersOnline) \textbf{5}, 86--89
  (2013).
\newblock \doi{10.3182/20130703-3-FR-4039.00028}

\bibitem{KiselevaKLN-2012-IEEE}
Kiseleva, M.A., Kuznetsov, N.V., Leonov, G.A., Neittaanm{\"{a}}ki, P.: Drilling
  systems failures and hidden oscillations.
\newblock In: IEEE 4th International Conference on Nonlinear Science and
  Complexity, NSC 2012 - Proceedings, pp. 109--112 (2012).
\newblock \doi{10.1109/NSC.2012.6304736}

\bibitem{kundur1994power}
Kundur, P., Balu, N.G., Lauby, M.G.: Power system stability and control,
  vol.~7.
\newblock McGraw-hill New York (1994)

\bibitem{Kurzin}
Kurzin, V., Seleznev, V.: Mechanism of emergence of intense vibrations of
  turbines on the {S}ayano-{S}hushensk hydro power plant.
\newblock Journal of Applied Mechanics and Technical Physics \textbf{51}(4),
  590--597 (2010)

\bibitem{KuznetsovL-2014-IFACWC}
Kuznetsov, N., Leonov, G.: Hidden attractors in dynamical systems: systems with
  no equilibria, multistability and coexisting attractors.
\newblock IFAC Proceedings Volumes (IFAC-PapersOnline) \textbf{19}, 5445--5454
  (2014).
\newblock \doi{10.3182/20140824-6-ZA-1003.02501}

\bibitem{KuznetsovKL-2013-DEDS}
Kuznetsov, N.V., Kuznetsova, O.A., Leonov, G.A.: Visualization of four normal
  size limit cycles in two-dimensional polynomial quadratic system.
\newblock Differential equations and dynamical systems \textbf{21}(1-2), 29--34
  (2013).
\newblock \doi{10.1007/s12591-012-0118-6}

\bibitem{KuznetsovLV-2010-IFAC}
Kuznetsov, N.V., Leonov, G.A., Vagaitsev, V.I.: Analytical-numerical method for
  attractor localization of generalized {C}hua's system.
\newblock IFAC Proceedings Volumes (IFAC-PapersOnline) \textbf{4}(1), 29--33
  (2010).
\newblock \doi{10.3182/20100826-3-TR-4016.00009}

\bibitem{kuznetsov2013elements}
Kuznetsov, Y.A.: Elements of applied bifurcation theory, vol. 112.
\newblock Springer Science \& Business Media (2013)

\bibitem{2001-Springer-Leonhard}
Leonhard, W.: Control of electrical drives.
\newblock Springer (2001)

\bibitem{LeonovBK-2010-DAN}
Leonov, G.A., Bragin, V.O., Kuznetsov, N.V.: Algorithm for constructing
  counterexamples to the {K}alman problem.
\newblock Doklady Mathematics \textbf{82}(1), 540--542 (2010).
\newblock \doi{10.1134/S1064562410040101}

\bibitem{leonov2009}
Leonov, G.A., Kondrat'eva, N.V.: Stability analysis of electric alternating
  current machines.
\newblock SPb: Isd. St. Petersburg. univ (2009)

\bibitem{LeonovK-2009}
Leonov, G.A., Kondrat'eva, N.V.: Stability analysis of electric alternating
  current machines [in {R}ussian].
\newblock SPb: Isd. St. Petersburg. univ (2009)

\bibitem{LeonovK-2011-DAN}
Leonov, G.A., Kuznetsov, N.V.: Algorithms for searching for hidden oscillations
  in the {A}izerman and {K}alman problems.
\newblock Doklady Mathematics \textbf{84}(1), 475--481 (2011).
\newblock \doi{10.1134/S1064562411040120}

\bibitem{LeonovK-2013-IJBC}
Leonov, G.A., Kuznetsov, N.V.: Hidden attractors in dynamical systems. {F}rom
  hidden oscillations in {H}ilbert-{K}olmogorov, {A}izerman, and {K}alman
  problems to hidden chaotic attractors in {C}hua circuits.
\newblock International Journal of Bifurcation and Chaos \textbf{23}(1) (2013).
\newblock \doi{10.1142/S0218127413300024}.
\newblock {a}rt. no. 1330002

\bibitem{LeonovKKSZ-2014}
Leonov, G.A., Kuznetsov, N.V., Kiseleva, M.A., Solovyeva, E.P., Zaretskiy,
  A.M.: Hidden oscillations in mathematical model of drilling system actuated
  by induction motor with a wound rotor.
\newblock Nonlinear Dynamics \textbf{77}(1-2), 277--288 (2014).
\newblock \doi{10.1007/s11071-014-1292-6}

\bibitem{LeonovKS-2015}
Leonov, G.A., Kuznetsov, N.V., Solovyeva: A simple dynamical model of
  hydropower plant: stability and oscillations.
\newblock Preprints of 1st IFAC Conference on Modelling, Identification and
  Control of Nonlinear Systems pp. 666--671 (2015)

\bibitem{LeonovKV-2011-PLA}
Leonov, G.A., Kuznetsov, N.V., Vagaitsev, V.I.: Localization of hidden {C}hua's
  attractors.
\newblock Physics Letters A \textbf{375}(23), 2230--2233 (2011).
\newblock \doi{10.1016/j.physleta.2011.04.037}

\bibitem{LeonovKV-2012-PhysD}
Leonov, G.A., Kuznetsov, N.V., Vagaitsev, V.I.: Hidden attractor in smooth
  {C}hua systems.
\newblock Physica D: Nonlinear Phenomena \textbf{241}(18), 1482--1486 (2012).
\newblock \doi{10.1016/j.physd.2012.05.016}

\bibitem{Merkurev}
Merkur'ev, G., Shargin, Y.: Stability of power systems.
\newblock {N}{O}{U} {C}{P}{K}{E}, Saint-Petersburg (2008)

\bibitem{nicolet2007}
Nicolet, C., Greiveldinger, B., H{\'e}rou, J.J., Kawkabani, B., Allenbach, P.,
  Simond, J.J., Avellan, F.: High-order modeling of hydraulic power plant in
  islanded power network.
\newblock Power Systems, IEEE Transactions on \textbf{22}(4), 1870--1880 (2007)

\bibitem{Pervozvanskii}
Pervozvanskii, A.: The automatic control theory course.
\newblock Nauka, Moscow (1986)

\bibitem{Act}
Rostechnadzor: The act of technical investigation into the causes of accident,
  occured on august 17, 2009 in the branch of the open joint-stock company
  "{R}us{H}ydro" - "{S}ayano-{S}hushenskaya {G}{E}{S} {P}.~{S}. {N}eporozneg".
\newblock Tech. rep. (2009)

\bibitem{sattoufsimulation2014}
Sattouf, M.: Simulation model of hydro power plant using matlab/simulink.
\newblock Journal of Engineering Research and Applications \textbf{4}(1),
  295--301 (2014)

\bibitem{TarbouriechGarcia_book11}
Tarbouriech, S., Garcia, G., {Gomes da Silva Jr.}, J., Queinnec, I.: Stability
  and Stabilization of Linear Systems with Saturating Actuators.
\newblock Springer-Verlag, London (2011)

\bibitem{TarbouriechTurner09}
Tarbouriech, S., Turner, M.: Anti-windup design: an overview of some recent
  advances and open problems.
\newblock IET Control Theory Appl. \textbf{3}(1), 1--19 (2009).
\newblock \urlprefix\url{{https://lra.le.ac.uk/handle/2381/4813}}

\bibitem{1931-Tricomi}
Tricomi, F.: Sur une equation differetielle de l'electrotechnique.
\newblock C.R. Acad. Sci. Paris. T. 193  (1931)

\bibitem{Wagner2011}
Wagner, H.J., J.Mathur: Introduction to Hydro Energy Systems: Basics,
  Technology and Operatio.
\newblock Springer Science \& Business Media (2011)

\bibitem{1958-Yanko-Trinitskii}
Yanko-Trinitskii, A.A.: New method for analysis of operation of synchronous
  motor for jump-like loads [in {R}ussian].
\newblock M.-L.: GEI (1958)

\bibitem{ZaccarianTeel_11book}
Zaccarian, L., Teel, A.: Modern Anti-windup Synthesis: {C}ontrol Augmentation
  for Actuator Saturation.
\newblock Princeton University Press (2011)

\end{thebibliography}

\end{document}